\documentclass[11pt]{article}
\date{}
\makeatletter
\@addtoreset{equation}{section}
\makeatother
\makeatother
\usepackage{amssymb,amsmath}
\usepackage{amsmath}
\usepackage{color}
\usepackage[colorlinks,linkcolor=blue,citecolor=blue]{hyperref}
\usepackage{cite}
\usepackage{graphicx}
\usepackage{times}

 \numberwithin{equation}{section}
 \newtheorem{thm}{Theorem}[section]
\newtheorem{lem}[thm]{Lemma}
\newtheorem{rem}[thm]{Remark}
\newtheorem{pro}[thm]{Proposition}
\newtheorem{cor}[thm]{Corollary}

\allowdisplaybreaks

\begin{document}
\title{The existence of ground state solutions for   critical H\'{e}non equations in $\mathbb{R}^N$\footnote{Supported by  NSFC(12301258,12271373,12171326).}}
\author{Cong Wang$^1$ \ \ \ \ \ Jiabao Su$^2$\footnote{Corresponding author. E-mail address: sujb@cnu.edu.cn.}\\ 
{\small  $^1$School of  Mathematics, Southwest Jiaotong University}\\ {\small Chengdu  610031, People's Republic of China}\\
{\small  $^2$School of Mathematical Sciences,
Capital Normal University}\\ {\small Beijing 100048, People's Republic of China}}
\maketitle
 \begin{abstract}
  In this paper we confirm that $2^*(\gamma)=\frac{2(N+\gamma)}{N-2}$ with $\gamma>0$ is exactly the critical exponent for the embedding from $H_r^1(\mathbb{R}^N)$ into $L^q(\mathbb{R}^N;|x|^\gamma)$($N\geqslant 3$) (see \cite{2007SWW-1,2007SWW-2}) and name it as the upper H\'enon-Sobolev critical exponent.
   Based on this fact we study the ground state solutions of critical H\'enon equations in $\mathbb{R}^N$ via the Nehari manifold methods and the great idea of Brezis-Nirenberg in \cite{1983BN}. We establish the existence of  the positive radial ground state solutions for the problem with one single upper H\'enon-Sobolev critical exponent. We also deal with the existence of the nonnegative radial ground state solutions for the problems with multiple critical exponents, including Hardy-Sobolev critical exponents or Sobolev critical exponents or the upper H\'{e}non-Sobolev critical exponents. 
 \end{abstract}
 \section{Introduction}
 In this paper we study the following equation
 \begin{eqnarray}
  \displaystyle \left\{\begin{array}{ll}
 \displaystyle -\Delta u+u=|x|^{\alpha}|u|^{2^*(\alpha)-2}u+\lambda|x|^{\beta}|u|^{p-2}u \ \ \ \hbox{in} \ \mathbb{R}^N, \\
 u\in H_r^1(\mathbb{R}^N),
  \end{array}
 \right.\label{1.1}
 \end{eqnarray}
  where $N \geqslant 3$,  $\alpha>0$,  $2^*(\alpha):= \frac{2(N+\alpha)}{N-2}$,
  $\lambda \in \mathbb{R}$ is a parameter, and $p$ satisfies
  \begin{eqnarray*}
  \left\{\begin{array}{ll}
 2_*(\beta)<p<2^*(\beta),    &  \beta\geqslant0,\\
 2_*(\beta)\leqslant p<2^*(\beta), \ & \ 0> \beta >-2,
  \end{array}
  \right.
   \  \ \ 2_*(\beta):=\left\{\begin{array}{ll}
 \frac{2(N-1+\beta)}{N-1},    &  \beta>0,\\
 2, \ & \ 0\geqslant \beta >-2.
  \end{array}
  \right.
  \end{eqnarray*}
  The Sobolev spaces of radial functions $H_r^1(\mathbb{R}^N)$ and $D_r^{1,2}(\mathbb{R}^N)$ are the completion of $C_{0,r}^\infty(\mathbb{R}^N)$ under the norms
  $$\|u\|:=\left(\int_{\mathbb{R}^N}(|\nabla u|^2+|u|^2) dx\right)^{\frac{1}{2}}, \ \ \ \|u\|_{D_r^{1,2}}:= \left(\int_{\mathbb{R}^N} |\nabla u|^2 dx\right)^\frac12. $$
  The well-known Sobolev embedding theorems tell us that the embedding from $H_r^1(\mathbb{R}^N)$ into $L^q(\mathbb{R}^N)$ is compact for all  $q\in \left(2, 2^*\right)$ (see \cite{1977Strauss,1996Willem}) and it is not true for $q = 2$
 and $q = 2^*$ where $2^*=\frac{2N}{N-2}$ is the classical Sobolev critical exponent.  We may refer $2^*$ and $2$ as the upper critical exponent and the lower critical exponent for the embedding $H_r^1(\mathbb{R}^N) \hookrightarrow L^q(\mathbb{R}^N)$. For $\gamma\geqslant-2$, we  define
      $$L^q(\mathbb{R}^N; |x|^\gamma):=\left\{ u: \mathbb{R}^N  \to \mathbb{R} \ \ \hbox{is Lebesgue measueable}, \ \int_{\mathbb{R}^N} |x|^\gamma |u|^q  dx <\infty  \right\}.$$By
      \cite[Theorem 1]{2007SWW-1,2007SWW-2} the embedding
   $H_r^1(\mathbb{R}^N) \hookrightarrow L^q(\mathbb{R}^N, |x|^\gamma)$
   is continuous for $2_*(\gamma) \leqslant q \leqslant 2^*(\gamma)$
   and it is compact for $2_*(\gamma) < q < 2^*(\gamma)$. Moreover, by \cite[Theorem 3.4]{2011SW}, we know that the embedding is also compact as $q=2_*(\gamma)=2$ with $-2<\gamma<0$, see Corollary \ref{cor1} in Section 2. It follows from these facts that
   \begin{eqnarray} \label{1.2}
   \Phi (u)= \frac{1}{2}\|u\|^2 - \frac{1}{2^*(\alpha)} \int_{\mathbb{R}^N} |x|^\alpha |u|^{2^*(\alpha)}dx - \frac{\lambda}{p} \int_{\mathbb{R}^N} |x|^\beta  |u|^p dx
 \end{eqnarray}
 is a well-defined $C^2$ functional on $H_r^1(\mathbb{R}^N)$ for $2_*(\beta)<p<2^*(\beta)$ with $\beta\geq0$ or $2\leqslant p<2^*(\beta)$ with $-2<\beta<0$. Thus the
 critical points of $\Phi$ are exactly the solutions of \eqref{1.1}.

 We confirm by Theorem \ref{thm3} in Section 2 that
 $2^*(\gamma)$ is exactly the upper critical exponent of the embedding from $H_r^1(\mathbb{R}^N)$  into $ L^q(\mathbb{R}^N; |x|^\gamma)$, this means that there is no embedding from  $H_r^1(\mathbb{R}^N)$  into $ L^q(\mathbb{R}^N; |x|^\gamma)$ for any $q>2^*(\gamma)$ and
 $H_r^1(\mathbb{R}^N) \hookrightarrow L^{2^*(\gamma)}(\mathbb{R}^N; |x|^\gamma)$ is not compact.

 The equation \eqref{1.1} is referred as a critical H\'enon equation on $\mathbb{R}^N$ since there is
 a critical term $|x|^\alpha |u|^{2^*(\alpha)-2} u$ with $\alpha>0$ contained in the equation.
 H\'enon equation is concerned with a semilinear equation on the unit ball $B=\{x\in \mathbb{R}^N: |x|<1\}$  with the weight $|x|^\alpha$ \begin{eqnarray}
  \left\{
  \begin{array}{ll}
    -\Delta u= |x|^\alpha u^{p-1}   & {\rm in} \ B,\\
 u>0 &{\rm in} \ B, \\
 u=0 &{\rm on} \ \partial B.
  \end{array}
 \right.\label{1.3}
 \end{eqnarray}
    The equation \eqref{1.3} with $\alpha>0$ was introduced by M. H\'enon in \cite{1973Henon} in studying the rotating stellar structures in 1973.  In the paper \cite{2002SSW}, Smets, Su and Willem first applied the variational methods to the H\'enon equation \eqref{1.3} and proved that there was $\alpha^*>0$ such that for $\alpha>\alpha^*$ and any $2<p<2^*$, the ground state solution of \eqref{1.3} was non-radial.  It follows that for large $\alpha>0$ and subcritical power $p$ the equation \eqref{1.3} has two   solutions in which one is radial and another one is non-radial. In \cite{1982Ni}, Ni proved that \eqref{1.3} had a radial solution for $2<p<2^*(\alpha)$.

 It is known in the literature that for $\gamma=-2$, $2^*(\gamma)=2$ is regarded as the Hardy critical exponent,  while for $-2<\gamma<0$,  $2^*(\gamma)$ is regarded as the Hardy-Sobolev critical exponent.
  For $\gamma>0$, we regard  $2^*(\gamma)$ as the upper H\'{e}non-Sobolev critical exponent of the embedding $H_r^1(\mathbb{R}^N) \hookrightarrow L^q(\mathbb{R}^N;|x|^\gamma)$ according to Theorem \ref{thm3} in Section 2. Thus $2^*(\gamma)=\frac{2(N+\gamma)}{N-2}$ is uniformly critical for all $\gamma\geqslant-2$.

In the celebrating paper \cite{1983BN}, Brezis and Nirenberg studied the following problem
 \begin{eqnarray}
 \left \{\begin{array}{ll}
 \displaystyle -\Delta u=u^{2^*-1}+ f(x, u) \ & \mbox{on} \ \Omega,\\
 \displaystyle u>0  \ &\mbox{on} \ \Omega,\\
 \displaystyle u=0 \ &\mbox{on} \ \partial \Omega
 \end{array}\label{1.4}
 \right.
 \end{eqnarray}
 with  $\Omega $ being a smooth bounded  domain in $\mathbb{R}^N (N\geqslant3)$.
 For $f(x, u)\equiv0$ and $\Omega$ is a star-shaped domain, it follows from Pohozaev's identity (see\cite{1965Pohozaev}) that the equation \eqref{1.4} has no solutions. For a general $f$, since
  the embedding from $H_0^1(\Omega)$ into $L^{2^*}(\Omega)$ is not compact,
  the energy functional associated to \eqref{1.4} may not satisfy the global Palais-Smale condition.
  With the help of lower-order perturbation and the mountain pass theorem\cite{1973AR},
  Brezis and Nirenberg \cite{1983BN} established the existence of a solution of \eqref{1.4} by constructing a minimax  value $c$ which was located in an interval $(-\infty, c^*)$ of the levels where the energy functional of \eqref{1.4} satisfied (PS) and $c^*>0$ was related to the best Sobolev constant.

    The great idea of Brezis and Nirenberg
    in \cite{1983BN} has motivated thousands of research works about the Brezis--Nirenberg problems
    in variant variational settings involving with different critical exponent problems. For the Brezis-Nirenberg problems on bounded domain with Hardy critical exponent $(\gamma=-2)$ or Hardy--Sobolev ($-2<\gamma<0$) critical exponent one can refer to the works \cite{2001CW,2010CL,1992Egnell,2001FG,2000GY,2004GK,2017GR2,2010HL,2011HL,2012Li-Lin} and the references therein.

   More recently, Wang and Su in \cite{2019Wang-Su} apply the ideas of Brezis and Nirenberg in \cite{1983BN} to deal with the  critical H\'enon equation on the unit ball
    \begin{eqnarray}
  \left\{
  \begin{array}{ll}
    -\Delta u= |x|^\alpha u^{2^*(\alpha)-1}+f(x, u)   & {\rm in} \ B,\\
 u>0 &{\rm in} \ B, \\
 u=0 &{\rm on} \ \partial B.
  \end{array}
 \right.\label{1.5}
 \end{eqnarray} In \cite{2019Wang-Su} the name of H\'enon-Sobolev critical exponent $2^*(\gamma)=\frac{2(N+\gamma)}{N-2}$ for $\gamma>0$ was given according to \cite[Theorem 2.2]{2019Wang-Su} which confirmed that $2^*(\gamma)$ is the critical exponent for the embedding $H_{0, r}^1(B) \hookrightarrow L^q(B, |x|^\gamma)$. In \cite{2019Wang-Su}, semilinear elliptic equations on the unit ball with multiple various critical exponents  have been studied.

  Motivated by \cite{1983BN} and \cite{2019Wang-Su}, in the present paper,  we consider the nontrivial radial  solutions for the critical H\'enon equation \eqref{1.1} on the whole spatial space $\mathbb{R}^N$. There is not nonzero solution  for \eqref{1.1} for $\lambda=0$ due to one reason that the embedding $H_r^1(\mathbb{R}^N) \hookrightarrow L^{2^*(\alpha)}(\mathbb{R}^N, |x|^\alpha)$ is not compact (Theorem \ref{thm3} in Section 2). Under the perturbation of subcritical term $|x|^\beta|u|^{p-2}u$ with $2_*(\beta)<p<2^*(\beta)$ and $\beta\geqslant0$ or $2_*(\beta)\leqslant p<2^*(\beta)$ and $-2<\beta<0$(see Corollary \ref{cor1} in Section 2), the functional $\Phi$ defined by \eqref{1.2} may satisfies (PS) at the levels below a positive number related to the best H\'enon--Sobolev constant $S_\alpha$ (see \cite{2013GGN,1983Lieb,2019Wang-Su}) for the embedding $D_r^{1,2}(\mathbb{R}^N)\hookrightarrow L^{2^*(\alpha)}(\mathbb{R}^N, |x|^\alpha)$.
  The existence of ground state solutions will be obtained by minimizing the energy functional $\Phi$ constrained on the Nehari manifold related to \eqref{1.1}.

 The paper is organized as follows. In Section 2, we recall by Theorem \ref{thm1} a weighted Sobolev embedding theorem for radial functions
 in $\mathbb{R}^N$ from \cite{2007SWW-1, 2007SWW-2}. As a special case of Theorem \ref{thm1} we obtain Corollary \ref{cor1} which states the   embedding from $H_r^1(\mathbb{R}^N)$ into $L^q(\mathbb{R}^N, |x|^\gamma)$ is compact for all $2_*(\gamma)<q<2^*(\gamma)$ with $\beta\geqslant0$ or $2_*(\gamma)\leqslant q<2^*(\gamma)$ with $-2<\beta<0$. Furthermore we confirm by Theorem \ref{thm3} that $2^*(\gamma)$ is exactly the (upper) critical exponent of embedding and is named to be the upper H\'{e}non-Sobolev critical exponent for $\gamma>0$. This result is quite new and ensures us to deal with elliptic equations with critical exponents of H\'enon-Sobolev type.
 In Section 3, we prove the existence of ground state solutions of problem with single H\'{e}non-Sobolev critical exponent using the Nehari manifold.  In Section 4 we deal with the existence of nontrivial ground state solutions for semilinear equations with double critical exponents, which may be Hardy-Sobolev critical exponent or Sobolev critical exponent or H\'{e}non-Sobolev critical exponent.

 \section{The H\'enon-Sobolev critical exponent. }
  Let us start this section with considering the radial solutions for the following equation
  \begin{eqnarray} \label{2.1}
  \left\{
  \begin{array}{ll}
  -\Delta u+u=Q(|x|)|u|^{q-2}u  & \mbox{in}\ \mathbb{R}^N,\\
   u(|x|)\rightarrow0  & \mbox{as} \ |x|\rightarrow\infty,
  \end{array}
  \right.
  \end{eqnarray}
 where the function $Q(r)>0$ is continuous on  $(0, \infty)$  and  satisfies the assumption
 \begin{enumerate}
   \item[$(\mathrm{Q})$] \ there exist real numbers $b_0$ and $b$ such that
 \begin{eqnarray*}
 \displaystyle \limsup_{r\rightarrow0}\frac{Q(r)}{r^{b_{0}}}<\infty, \ \ \ \ \limsup_{r\rightarrow\infty}\frac{Q(r)}{r^{b}}<\infty.
 \end{eqnarray*}
 \end{enumerate}
  Define for $q\geqslant 1$,
  $$L^q(\mathbb{R}^N; Q(|x|)):= \left\{u:\mathbb{R}^N \to \mathbb{R} \ \  \mbox{is measurable,} \ \int_{\mathbb{R}^N}
  Q(|x|)|u|^q dx<\infty\right\}.$$
  Define \begin{eqnarray} \label{2.2}
2_*(b):=\left\{
  \begin{array}{ll}
\frac{2(N-1+b)}{N-1} \ &{\rm for}   \   b>0, \\
2\ & {\rm for } \ -2\leqslant b \leqslant 0,
  \end{array} \ \
\right.2^*(b_0):=\frac{2(N+b_0)}{N-2} \  \ \hbox{for} \  b_0 \geqslant -2.
\end{eqnarray}
 As a special case of Su, Wang and Willem \cite[Theorem 1]{2007SWW-1,2007SWW-2}, Su and Wang \cite[Theorem 3.4]{2011SW}, we have the following embedding result.
 \begin{thm}\label{thm1}
 Assume $(\mathrm{Q})$ holds with $b\geqslant-2$ and $b_0\geqslant-2$ be such that $2_*(b) < 2^*(b_0)$. Then the embedding
 $$H_r^1(\mathbb{R}^N)\hookrightarrow L^q(\mathbb{R}^N;Q)$$ is continuous
  for $2_*(b) \leqslant  q \leqslant  2^*(b_0)$. Furthermore, the embedding is compact for $2_*(b)<q<2^*(b_0)$ and compactness still holds as $p=2_*(b)=2$ with $-2<b<0$.
 \end{thm}
 In the case that $Q(|x|)= |x|^{\gamma}, b_0=b=\gamma$,   we   get the following corollary.
 \begin{cor}\label{cor1} Assume that $\gamma \geqslant -2$. The embedding
 \begin{eqnarray} \label{2.3}
 H_r^1(\mathbb{R}^N) \hookrightarrow L^q(\mathbb{R}^N;|x|^\gamma)
 \end{eqnarray}
 is continuous for $2_*(\gamma) \leqslant  q \leqslant  2^*(\gamma)$, and it is compact for $2_*(\gamma)<q<2^*(\gamma)$ and compactness is true as $q=2_*(\gamma)=2$ with $-2<\gamma<0$.
\end{cor}
Applying the above conclusion, we can define well the "first" eigenvalue:
\begin{eqnarray*}
\lambda_{1\gamma}:=\inf_{u\in H_r^1(\mathbb{R}^N)\setminus\{0\}}\frac{\int_{\mathbb{R}^N}|\nabla u|^2+|u|^2dx}{\int_{\mathbb{R}^N}|x|^\gamma |u|^2dx}\ {\rm with}\ -2<\gamma<0
\end{eqnarray*}
and it can be achieved by positive eigenfunction $\varphi_{1\gamma}$ combining with strongly maximum principle.

 Furthermore, by scaling arguments, we can prove that $2^*(\gamma)$ is the upper critical exponent for the embedding \eqref{2.3}. That is the following conclusion.
\begin{thm}\label{thm3}
Assume that $\gamma \geqslant -2$.  For any $q>2^*(\gamma)$, there is no embedding from $H_r^1(\mathbb{R}^N)$ into $L^q(\mathbb{R}^N; |x|^\gamma)$. The embedding $H_r^1(\mathbb{R}^N) \hookrightarrow L^{2^*(\gamma)}(\mathbb{R}^N;|x|^\gamma)$ is not compact.
\end{thm}
{\bf Proof.}  \  We only need to construct a counter-example to illustrate the conclusions of the theorem.
 For $k\in \mathbb{N}$, we define a sequence of radial functions $\{u_k\}_{k=1}^\infty$ as follows.
 $$u_k(|x|)=k^{\frac{N-2}{4}}e^{-\frac{k|x|^2}{2\cdot2^*(\gamma)}}, \ \ k\in \mathbb{N}. $$
 Direct computation shows that
 \begin{eqnarray*}
 \int_{\mathbb{R}^N}|\nabla u_k|^2dx=\frac{(2^*(\gamma))^{\frac{N-2}{2}}\omega_N}{2}\Gamma\left(\frac{N+2}{2}\right),
 \end{eqnarray*} \begin{eqnarray*}
\int_{\mathbb{R}^N}|u_k|^2 dx \leqslant \frac{(2^*(\gamma))^{\frac{N}{2}}\omega_N }{2}\Gamma\left(\frac{N}{2}\right),
\end{eqnarray*} where $\omega_N$ is surface area of the unite sphere in $\mathbb{R}^{N}$ and $\Gamma$  is the gamma function.
 Therefore $\{u_k\} \subset H_r^1(\mathbb{R}^N)$ is a bounded sequence.

For $q>2^*(\gamma)$, we have
\begin{eqnarray*}
\int_{\mathbb{R}^N}|x|^{\gamma}|u_k|^q dx =k^{\frac{q(N-2)}{4}-\frac{N+\gamma}{2}}\frac{\omega_N }{2}\left(\frac{2\cdot2^*(\gamma)}{q}\right)
^{\frac{N+\gamma}{2}}\Gamma\left(\frac{N+\gamma}{2}\right).
\end{eqnarray*}
 Since $q>2^*(\gamma)$ is equivalent to $\frac{q(N-2)}{4}-\frac{N+\gamma}{2}>0$, it follows that
\begin{eqnarray*}
\int_{\mathbb{R}^N}|x|^\alpha|u_k|^q dx \rightarrow   \infty \ \  \hbox{as} \ k \to\infty, \ \ \ \forall \ q>2^*(\gamma).
\end{eqnarray*}
Therefore there is no embedding from $H_r^1(\mathbb{R}^N)$ into  $L^q(\mathbb{R}^N;|x|^{\gamma})$ for all
 $q>2^*(\gamma)$.

On one hand, we have
\begin{eqnarray*}
\int_{\mathbb{R}^N}|x|^{\gamma}|u_k|^{2^*(\gamma)} dx =2^{\frac{N+\gamma-2}{2}}\omega_N \Gamma\left(\frac{N+\gamma}{2}\right)>0.
\end{eqnarray*}
On the other hand, it is easy to see that
$$u_k(|x|)\rightarrow 0 \ \ \ \mbox{for \ a.e. } \  x \in \mathbb{R}^N, \ \ \ k\rightarrow\infty.$$
Therefore $\{u_k\}$ does not contain any subsequence converging in $L^{2^*(\gamma)}(\mathbb{R}^N, |x|^{\gamma})$. $\hfill \Box$

For $\gamma>0$ we name $2^*(\gamma)$ as the  upper H\'enon-Sobolev critical exponent for the embedding from $H_r^1(\mathbb{R}^N)$ into $L^q(\mathbb{R}^N; |x|^\gamma)$. It is an open problem whether or not $2_*(\gamma)=\frac{2(N-1+\gamma)}{N-1}$ is the
lower critical exponent for this embedding.

By Theorem \ref{thm3}, since the embedding $H_r^1(\mathbb{R}^N)\hookrightarrow L^{2^*(\alpha)}(\mathbb{R}^N;|x|^\alpha)$ is not compact with $\alpha>0$, it follows that the
energy functional $\Phi$ of \eqref{1.1} may not satisfy the global Palais-Smale condition. However, it may satisfy the ${\rm(PS)_c}$ condition at the energy levels $c$ in certain intervals. It concerns with     the following equation
\begin{eqnarray}
\left\{\begin{array}{ll} \label{2.4}
-\Delta u=|x|^\alpha u^{2^*(\alpha)-1} & \text{in} \ \mathbb{R}^N, \\
 u>0 & \text {in} \ \mathbb{R}^{N}, \\
u\in D_r^{1,2}(\mathbb{R}^{N}).
 \end{array}
 \right.
\end{eqnarray}
By \cite{1981GS, 2013GGN, 1983Lieb}, \eqref{2.4} has a unique(up to dilations) radial solution given by
\begin{eqnarray}
U_{\epsilon, \alpha}(x)=\frac{C(\alpha,N)\epsilon^{\frac{N-2}{2}}}{\left(\epsilon^{2+\alpha}+|x|^{2+\alpha}\right)
^{\frac{N-2}{2+\alpha}}},\label{2.5}
\end{eqnarray}
where $C(\alpha,N)=[(N+\alpha)(N-2)]^{\frac{N-2}{4+2\alpha}}$.  They are extremal functions for the following inequality
\begin{eqnarray}
\int_{\mathbb{R}^{N}}|\nabla u|^{2} dx \geqslant S_\alpha\left(\int_{\mathbb{R}^{N}}|x|^{\alpha}|u|^{2^*(\alpha) } dx \right)^{\frac{2}{2^*(\alpha)}},
\ u\in D^{1,2}_r(\mathbb{R}^{N}).\label{2.6}
\end{eqnarray}
 We call \eqref{2.6}  the ``H\'{e}non-Sobolev'' inequality in which $S_\alpha$ can be written as
 \begin{eqnarray*}
 S_\alpha:=S_\alpha(\mathbb{R}^N)=(N+\alpha)(N-2)\left(\frac{\omega_N}{2+\alpha}
 \cdot \frac{\Gamma^2\left(\frac{N+\alpha}{2+\alpha}\right)}
 {\Gamma\left(\frac{2(N+\alpha)}{2+\alpha}\right)}\right)^{\frac{2+\alpha}{N+\alpha}}.
 \end{eqnarray*}
 See \cite{2019Wang-Su} for a computation. When $\alpha=0$, $S_\alpha$ coincides the best Sobolev constant $S_0$(see \cite{1976Talenti}) for $\alpha=0$, which is
\begin{eqnarray*}
S_0=(N)(N-2)\left(\frac{\omega_N}{2}
\cdot \frac{\Gamma^2\left(\frac{N}{2}\right)}
{\Gamma\left(N\right)}\right)^{\frac{2}{N}}=\pi N(N-2)\left(\frac{\Gamma\left(\frac{N}{2}\right)}
{\Gamma\left(N\right)}\right)^{\frac{2}{N}}.
\end{eqnarray*} and when $-2<\alpha<0$,
 then $S_\alpha(\mathbb{R}^N)$ is the best Hardy-Sobolev constant(see \cite{2016GR1}).

 We can construct a smooth function on a given ball $B_r=\{x\in \mathbb{R}^N:  |x|< r\}$ from the function \eqref{2.5} for  later uses.  Define
 \begin{eqnarray} \label{2.7}
 u_{\epsilon,\alpha}= \varphi(x)U_{\epsilon, \alpha}(x),
 \end{eqnarray}
 where $\varphi(|x|)\in C_{0,r}^{\infty} (B_r)$,
 $0 \leqslant \varphi(|x|) \leqslant 1$ and $\varphi(|x|)$ satisfies
 \begin{eqnarray*}
 \left\{\begin{array}{ll}
 \varphi(|x|)\equiv 1 &  \mbox{for} \ |x| \leqslant  R,\\
\varphi(|x|)\equiv0 &  \mbox{for} \ |x| \geqslant 2R.
 \end{array}
 \right.
\end{eqnarray*}
By careful computations, we have
\begin{eqnarray}
&&\int_{\mathbb{R}^N}\left|\nabla u_{\epsilon,\alpha}\right|^2dx
=S_\alpha^{\frac{N+\alpha}{2+\alpha}}+O\left(\epsilon^{N-2}\right). \label{2.8}\\
&&\int_{\mathbb{R}^N}|x|^\alpha |u_{\epsilon,\alpha}|^{2^*(\alpha)}dx
=S_\alpha^{\frac{N+\alpha}{2+\alpha}}+O\left(\epsilon^{N+\alpha}\right).\label{2.9}
\end{eqnarray}
\begin{eqnarray}
\int_{\mathbb{R}^N}|u_{\epsilon,\alpha}|^2dx=\left \{\begin{array}{ll}
O(\epsilon), &  N=3, \\
O(\epsilon^2)+C\epsilon^2|\ln \epsilon|, &  N=4,\\
 O(\epsilon^{N-2})+C\epsilon^2, &  N\geqslant 5.
\end{array}
\right.\label{2.10}
\end{eqnarray}
\begin{eqnarray}
\int_{\mathbb{R}^N}|x|^\beta|u_{\epsilon,\alpha}|^pdx=\left \{\begin{array}{ll}
 O(\epsilon^{\frac{N+\beta}{2}})+C\epsilon^{\frac{N+\beta}{2}}|\ln \epsilon|, &   p=\frac{N+\beta}{N-2},\\
 O(\epsilon^{\frac{p(N-2)}{2}})+C\epsilon^{N+\beta-\frac{p(N-2)}{2}}, &   p>\frac{N+\beta}{N-2},
\end{array}
\right.\label{2.11}
\end{eqnarray}
where $\beta\geqslant-2$.
The function $u_{\epsilon,\alpha}$  and the estimates \eqref{2.8}--\eqref{2.11} will be used   in the below.


\section{Problems with single H\'{e}non-Sobolev critical exponent}

In this section we study the existence of ground state solutions of the equation
  \begin{eqnarray}
\left \{\begin{array}{ll}
\displaystyle -\Delta u+u=|x|^{\alpha}|u|^{2^*(\alpha)-2}u+\lambda |x|^\beta |u|^{p-2}u  & \hbox{in}  \ \  {\mathbb{R}^N}, \\
\displaystyle u\in H^1_r({\mathbb{R}^N}), &
\end{array}\label{3.1}
\right.
\end{eqnarray}
 where $N\geqslant 3$, $\alpha>0$, $\beta>-2$, $2_*(\beta) < p < 2^*(\beta)$ with $\beta\geq0$ or $2_*(\beta)\leqslant p < 2^*(\beta)$ with $-2<\beta<0$, and $2^*(\star)$ is the upper H\'{e}non-Sobolev critical exponent.
 It was proved in \cite{1983B-Lions} via the Pohozaev's identity(\cite{1965Pohozaev}) that the equation
 \begin{eqnarray}
 \left\{
  \begin{array}{ll}
  -\Delta u+u=|u|^{2^*-2} u,   &  \hbox{in} \ \mathbb{R}^N, \\
 u\in H^1(\mathbb{R}^N), \ \ u \not=0 &
  \end{array}
 \right.\label{3.2}
 \end{eqnarray}
 has no solutions.  We consider the following equation with  $\alpha>0$
 \begin{eqnarray}
\left \{\begin{array}{ll}
\displaystyle -\Delta u+u=|x|^{\alpha}|u|^{2^*(\alpha)-2}u,  &  \hbox{in} \  {\mathbb{R}^N},\\
\displaystyle u\in H^1_r({\mathbb{R}^N}), \  u\ne0. &
\end{array}\label{3.3}
\right.
\end{eqnarray}
We will have a same conclusion that \eqref{3.3} has no solutions. Assume $u\in H^1_r(\mathbb{R}^N)$ is a solution of \eqref{3.3}. Denote $a(|x|):=|x|^\alpha |u|^{2^*(\alpha)-2}-1$. For any ball $B_R\subset \mathbb{R}^N$, we have by Corollary \ref{cor1} that
\begin{eqnarray*}
\int_{B_R}|a(|x|)|^{\frac{N}{2}}dx& \leqslant & C \int_{B_R} |x|^{\frac{N\alpha}{2}} |u|^{\frac{2(N+\frac{N\alpha}{2})}{N-2}}dx+C|B_R|< \infty.
\end{eqnarray*}
Thus  $a \in L^{N/2}(B_R)$. According to the ideas of  Lemma B.3 in \cite{2008Struwe},
 we have $u \in L_r^q(B_R)$ for any $q<\infty$. For any $t<\infty$, we have
\begin{eqnarray*}
\int_{B_R}||x|^{\alpha}|u|^{2^*(\alpha)-2}u-u|^tdx &\leqslant  &C\int_{B_R}|x|^{t\alpha} |u|^{\frac{t(N+2\alpha+2)}{N-2}}dx+C\int_{B_R} |u|^tdx\\
& \leqslant & C R^{t \alpha} \int_{B_R} |u|^{\frac{t(N+2\alpha+2)}{N-2}}dx+C\int_{B_R} |u|^tdx\\
&<& \infty.
 \end{eqnarray*}
 Hence by Calderon-Zygmund inequality, $u\in W_r^{2,t}(B_R)$ for any $t<\infty$.
 It follows from Sobolev embedding theorem that $u\in C^{1, s}_{r}(B_R)$ for $0 \leqslant s<1$.
 Using the arguments in \cite{1983B-Lions}, the Pohozaev's identity of \eqref{3.3} reads as
 \begin{eqnarray}
 \int_{\mathbb{R}^N}|\nabla u|^2dx+\frac{N}{N-2}\int_{\mathbb{R}^N} |u|^2 dx = \int_{\mathbb{R}^N} |x|^{\alpha} |u|^{2^*(\alpha)} dx.\label{3.4}
\end{eqnarray} Since a solution $u \in H_r^1(\mathbb{R}^N)$ of \eqref{3.3} verifies
 \begin{eqnarray*}
 \|u\|^2=\int_{\mathbb{R}^N}(|\nabla u|^2+|u|^2)dx=\int_{\mathbb{R}^N}|x|^\alpha |u|^{2^*(\alpha)}dx,
 \end{eqnarray*} it follows that $\int_{\mathbb{R}^N} |u|^2dx=0$ and then $u=0$. A contradiction.

 It is easy to see that when   $2_*(\beta) < p < 2^*(\beta)$, \eqref{3.1} has no nontrivial solutions for $\lambda<0$. We only consider the equation \eqref{3.1} with the case $\lambda>0$.  We assume
 \begin{equation} \label{3.5}
  \left\{
  \begin{array}{ll}
 2(2+\beta)<p < 2^*(\beta), &  \ \beta \geqslant -1, \  N=3;\\
2<p<2^*(\beta), &  \ \beta<-1, \  N=3;\\
 2 < p < 2^*(\beta),& -2<\beta\leqslant 0, \ N \geqslant 4;\\
 \frac{2(N-2+\beta)}{N-2}< p < 2^*(\beta), & \beta>0, \ N\geqslant 4.
   \end{array}
\right.
\end{equation}
\begin{equation} \label{new3.6}
  \left\{
  \begin{array}{ll}
p=2, &  \ -2<\beta \leqslant -1, \  N=3;\\
p=2, &  \-2< \beta<0, \  N\geqslant4.
   \end{array}
\right.
\end{equation}
The main result in this section is the following theorem.
\begin{thm}\label{thm5} \
{\rm (i)} Assume $\beta>-2$ and $2_*(\beta)<p<2^*(\beta)$, then there exists $\lambda^*>0$ such that \eqref{3.1} has a positive ground state solution for $\lambda>\lambda^*$.

{\rm (ii)} \ Assume that $p$ satisfies \eqref{3.5}.  Then \eqref{3.1} has a positive ground state solution for any $\lambda>0$.

{\rm (iii)} \ Assume that $p$ satisfies \eqref{new3.6}.  Then \eqref{3.1} has a positive ground state solution for any $0<\lambda<\lambda_{1\beta}$.
\end{thm}
\begin{rem}\label{rem3.2}
We remark that for $N\geq3, -2<\beta<0, p=2$, the problem \eqref{3.1} has no positive solutions as $\lambda\geq\lambda_{1\beta}$. Indeed, let $\lambda\geq\lambda_{1\beta}$  and $u_0$ be a positive solution, then for the "first" positive eigenfunction $\varphi_{1\beta}$(see section {\rm 2}), we have
\begin{eqnarray*}
\lambda_{1\beta}\int_{\mathbb{R}^N} |x|^\beta \varphi_{1\beta} u_0dx&=&\int_{\mathbb{R}^N}\nabla u_0\nabla \varphi_{1\beta}+u_0\varphi_{1\beta}dx\\
&=&\int_{\mathbb{R}^N}|x|^\alpha u_0^{2^*(\alpha)-1}\varphi_{1\beta}dx+\lambda\int_{\mathbb{R}^N} |x|^\beta \varphi_{1\beta} u_0dx\\
&>&\lambda\int_{\mathbb{R}^N} |x|^\beta \varphi_{1\beta} u_0dx.
\end{eqnarray*}
It follows that $\lambda<\lambda_{1\beta}$, a contradiction leads to the problem \eqref{3.1} has no positive solutions.
\end{rem}
We  will apply the  Nehari manifold methods to prove this theorem.
The energy functional corresponding  to \eqref{3.1} is defined as
\begin{eqnarray}\label{3.6}
\Phi(u)=\frac{1}{2}\|u\|^2-\frac{1}{2^*(\alpha)}\int_{\mathbb{R}^N}|x|^\alpha |u|^{2^*(\alpha)}dx-\frac{\lambda}{p}
\int_{\mathbb{R}^N}|x|^\beta|u|^pdx, \ \ u\in H^1_r(\mathbb{R}^N).
\end{eqnarray}
By Theorem \ref{thm1} and Corollary \ref{cor1}, we have that $\Phi\in C^2(H_r^1(\mathbb{R}^N), \mathbb{R})$. The Nehari manifold of $\Phi$ is defined as
\begin{eqnarray*}
 \mathcal{N}=\left\{u\in H_r^1(\mathbb{R}^N)\setminus\{0\} \ | \ \Psi(u):= \langle \Phi'(u), u\rangle =0\right\}.
\end{eqnarray*}
 Define
 \begin{eqnarray} \label{3.7}
  \displaystyle m:=\inf_{u\in\mathcal{N}} \Phi (u).
 \end{eqnarray}
For the sake of conciseness, we give an assumption
\begin{eqnarray}
\left\{
  \begin{array}{ll}
2_*(\beta)<p<2^*(\beta)\ {\rm with}\ \beta>-2\ {\rm and}\  \lambda>0;\\
p=2\ {\rm with}\ -2<\beta<0\ {\rm and}\  0<\lambda<\lambda_{1\beta}.
  \end{array}
\right.\label{3.8}
\end{eqnarray}
  \begin{pro}\label{pro1}
Assume  \eqref{3.8} holds.
  \begin{enumerate}
 \item[{\rm (i) }] For $u\in H_r^1(\mathbb{R}^N)\backslash\{0\}$, there exists a unique $t_u>0$ such that  $t_uu\in\mathcal{N}$ and $\Phi (t_u u)=\max_{t\geqslant  0}\Phi (tu)$, moreover, the manifold $\mathcal{N}$ is nonempty;
     \item[{\rm (ii) }] The manifold $\mathcal{N}$ is $C^1$ regular, that is $\langle \Psi'(u), u\rangle=\langle \Phi''(u)u, u\rangle\ne0$ for any $u\in \mathcal{N}$;
     \item[{\rm (iii) }] $\mathcal{N}$ is closed and bounded away from $0$, and $m>0$.
 \end{enumerate}
  \end{pro}
{\bf Proof.} \ (i) \ For any $u\in H^1_r(\mathbb{R}^N)\backslash\{0\}$,
\begin{eqnarray}
\Phi (tu)=\frac{t^2}{2}\|u\|^2 -\frac{t^{2^*(\alpha)}}{2^*(\alpha)}\int_{\mathbb{R}^N}|x|^\alpha |u|^{2^*(\alpha)}dx-\frac{\lambda t^p}{p} \int_{\mathbb{R}^N}|x|^\beta|u|^pdx.\label{3.9}
\end{eqnarray}
For the case $p>2_*(\beta), \lambda>0$, since $\hat p= \min\{p, 2^*(\alpha)\}>2$, it follows that there is
  a unique $t_u>0$ such that
  \begin{eqnarray}
\Phi (t_u u)=\max_{t\geqslant 0} \Phi(tu), \ \ \ \frac{d\Phi(t u)}{dt}\bigg|_{t=t_u}= \langle \Phi'(t_u u), u\rangle = 0.\label{3.10}
  \end{eqnarray}
For the case $p=2$ and $0<\lambda<\lambda_{1\beta}$, then
\begin{eqnarray}
\int_{\mathbb{R}^N}|\nabla u|^2+|u|^2dx-\lambda\int_{\mathbb{R}^N}|x|^\beta|u|^2dx\geq
\frac{\lambda_{1\beta}-\lambda}{\lambda_{1\beta}}\int_{\mathbb{R}^N}|\nabla u|^2+|u|^2dx>0.\label{3.11}
\end{eqnarray}
Combining \eqref{3.9} with $p=2$  and \eqref{3.11}, we get that the fact \eqref{3.10} is valid. Thus  $t_u u\in \mathcal{N}$ and $\mathcal{N}$ is nonempty.

(ii) \ For   $u\in \mathcal{N}$, we have
\begin{eqnarray*}
\langle \Phi'(u), u\rangle = \|u\|^2 -\int_{\mathbb{R}^N}|x|^\alpha |u|^{2^*(\alpha)}dx-\lambda\int_{\mathbb{R}^N}|x|^\beta|u|^pdx=0.
\end{eqnarray*}
For $p>2_*(\beta)$ with $\lambda>0$, since $2< {\hat{p}}:=\min\{p,2^*(\alpha)\}$,
\begin{eqnarray}
\left.
  \begin{array}{ll}
\langle \Psi'(u), u\rangle &=\langle \Phi''(u)u, u\rangle-\hat{p}\langle \Phi'(u), u\rangle\\
&=(2-\hat{p}) \|u\|^2+(\hat{p}-2^*(\alpha))\int_{\mathbb{R}^N}|x|^\alpha |u|^{2^*(\alpha)}dx  +\lambda(\hat{p}-p) \int_{\mathbb{R}^N}|x|^\beta|u|^pdx\\
& \leqslant  (2-\hat{p}) \|u\|^2 <0.
  \end{array}
\right.\label{3.12}
\end{eqnarray}
For $p=2$ with $0<\lambda<\lambda_{1\beta}$, using the fact \eqref{3.11} and inequality \eqref{3.12} with $\hat{p}=2^*(\alpha),p=2$, then we can deduce that
\begin{eqnarray*}
\langle \Psi'(u), u\rangle\leq \frac{(2-2^*(\alpha))(\lambda_{1\beta}-\lambda)}{\lambda_{1\beta}}\int_{\mathbb{R}^N}|\nabla u|^2+|u|^2dx<0.
\end{eqnarray*}
It follows from the implicit function theorem that $\mathcal{N}$ is a $C^1$-manifold and is regular.

(iii) \ Let $u\in \mathcal{N}$. For the case $p>2_*(\beta)$, by the H\'enon-Sobolev inequality \eqref{2.6} and Theorem \ref{thm1}, one has
\begin{eqnarray}
\left.
  \begin{array}{ll}
0&=\|u\|^2-\int_{\mathbb{R}^N}|x|^\alpha |u|^{2^*(\alpha)}dx-\lambda\int_{\mathbb{R}^N}|x|^\beta|u|^pdx\\
&\geqslant \|u\|^2-S_\alpha^{-\frac{2^*(\alpha)}{2}}\|u\|^{2^*(\alpha)}-\lambda S_{p\beta }^p\|u\|^p.
  \end{array}
\right.\label{3.13}
\end{eqnarray}
where $S_{p\beta}>0$ is the embedding constant from $H_r^1(\mathbb{R}^N)$ into $L^p(\mathbb{R}^N, |x|^\beta)$. Hence there exists $\delta_1:=\delta_1(N, \alpha,\beta)>0$ such that
 $\|u\|\geqslant \delta_1$ for all $u\in \mathcal{N}.$
It follows that
\begin{eqnarray*}
\Phi(u) \geqslant \left(\frac{1}{2}-\frac{1}{\hat{p}}\right) \delta_1^2,  \ \ \ \ \ \forall \ u\in \mathcal{N}.
\end{eqnarray*}
For the case of $p=2$ and $0<\lambda<\lambda_{1\beta}$, similarly, combining with the facts \eqref{3.11}, \eqref{3.13} with $p=2$, we can find a $\delta_2>0$ such that
\begin{eqnarray*}
\Phi(u) \geqslant \left(\frac{1}{2}-\frac{1}{2^*(\alpha)}\right)\frac{\lambda-\lambda_{1\beta}}{\lambda_{1\beta}}\delta_2^2,  \ \ \ \ \ \forall \ u\in \mathcal{N}.
\end{eqnarray*}
Thus $m>0$. \hfill$\Box$

Suppose \eqref{3.8}, following the arguments of \cite[Chapter 4]{1996Willem}, we can prove that
\begin{eqnarray}  \label{3.14}
m= \hat c:=\inf_{u\in H_r^1(\mathbb{R}^N) \setminus\{0\}}\max_{t\geqslant 0}\Phi (t u).
\end{eqnarray}

\begin{lem}\label{lem1} Let condition \eqref{3.8} hold. If $\{u_n\}\subset H^1_r(\mathbb{R}^N)$ is a  ${{\rm(PS)}_c}$ sequence of $\Phi$ with $c<\frac{2+\alpha}{2(N+\alpha)}S_\alpha^{\frac{N+\alpha}{2+\alpha}}$, then it contains   a convergent subsequence.
\end{lem}
{\bf Proof.}  Let $\{u_n\}\subset H^1_r(\mathbb{R}^N)$ be such that
\begin{eqnarray*}
\Phi(u_n)\rightarrow c, \ \  \ \ \ \Phi'(u_n) \rightarrow 0 \ \ \mbox{as} \ \ n\rightarrow\infty.
\end{eqnarray*}
Then as $n\to\infty$
 \begin{eqnarray*}
 &&\Phi(u_n)=\frac{1}{2}\|u_n\|^2 -\frac{1}{2^*(\alpha)}\int_{\mathbb{R}^N}|x|^\alpha|u_n|^{2^*(\alpha)}dx
-\frac{\lambda}{p}\int_{\mathbb{R}^N}|x|^\beta|u_n|^pdx \to c,\\
&&\langle \Phi'(u_n), u_n\rangle=\|u_n\|^2 -\int_{\mathbb{R}^N}|x|^\alpha|u_n|^{2^*(\alpha)}dx
-\lambda\int_{\mathbb{R}^N}|x|^\beta|u_n|^pdx=o(1)\|u_n\|.
 \end{eqnarray*}
As $p>2_*(\beta)$, for $n \in \mathbb{N}$ large, taking $\hat{p}>2$ in Proposition \ref{pro1}, we have
\begin{eqnarray}
 \left.
  \begin{array}{ll}
& c+1+o(1)\|u_n\|\\[2mm]
\geqslant  & \displaystyle \Phi(u_n)-\frac{1}{\hat{p}}\langle \Phi'(u_n), u_n\rangle\\[2mm]
= & \displaystyle \left(\frac{1}{2}-\frac{1}{\hat{p}}\right)\|u_n\|^2+\left(\frac{1}{\hat{p}}
-\frac{1}{2^*(\alpha)}\right)\int_{{\mathbb{R}^N}}|x|^\alpha|u_n|^{2^*(\alpha)}dx
+\lambda\left(\frac{1}{\hat{p}}-\frac{1}{p}\right)\int_{{\mathbb{R}^N}}|x|^\beta|u_n|^pdx\\[2mm]
\geqslant & \displaystyle  \left(\frac{1}{2}-\frac{1}{\hat{p}}\right)\|u_n\|^2.
  \end{array}
\right.\label{3.15}
\end{eqnarray}
As $p=2_*(\beta)=2$ and $\lambda<\lambda_{1\beta}$. Using inequality \eqref{3.11} and \eqref{3.15} with $\hat{p}=2^*(\alpha), p=2$, we get
\begin{eqnarray}
 c+1+o(1)\|u_n\|\geqslant  \left(\frac{1}{2}
-\frac{1}{2^*(\alpha)}\right)\frac{\lambda_{1\beta}-\lambda}{\lambda_{1\beta}}\|u_n\|^2.\label{3.16}
\end{eqnarray}
It follows from \eqref{3.15} with $\hat{p}>2$ and \eqref{3.16} that $\{\|u_n\|\}$ is bounded. Going if necessary to a subsequence, we can assume that there is $u\in H_r^1(\mathbb{R}^N)$ such that
\begin{eqnarray}
  \left\{\begin{array}{ll}
&u_n\rightharpoonup u \ \ \mbox{in} \ H^1_r(\mathbb{R}^N);\\
&u_n\rightharpoonup u \ \ \mbox{in} \ L^{2^*(\alpha)}(\mathbb{R}^N;|x|^\alpha);\\
&u_n\rightarrow u \ \  \mbox{in} \ L^p(\mathbb{R}^N;|x|^\beta),  \ \ \mbox{for} \ 2_*(\beta)<p<2^*(\beta) \ \hbox{with} \ \beta>-2;\\
&u_n\rightarrow u \ \  \mbox{in} \ L^2(\mathbb{R}^N;|x|^\beta),  \ \ \mbox{for} \ p=2_*(\beta)=2 \ \hbox{with} \ -2<\beta<0;\\
&u_n(x)\rightarrow u(x) \ \mbox{a.e.\ on} \  \mathbb{R}^N.
  \end{array}
\right.\label{3.17}
\end{eqnarray}
By \eqref{3.17}  we have that
$$\int_{\mathbb{R}^N} |x|^\alpha|u_n|^{2^*(\alpha)-2}u_nwdx \rightarrow \int_{\mathbb{R}^N}|x|^\alpha |u|^{2^*(\alpha)-2}uwdx,   \ \ \ \forall \ w \in H^1_r(\mathbb{R}^N).$$
Therefore $u$ solves weakly the equation
\begin{eqnarray} \label{3.18}
 -\Delta u+u=|x|^{\alpha}|u|^{2^*(\alpha)-2}u+\lambda |x|^\beta u^p.
 \end{eqnarray}
Thus
\begin{eqnarray}
\Phi(u)&=&\frac{1}{2}\|u\|^2-\frac{1}{2^*(\alpha)}\int_{\mathbb{R}^N}|x|^\alpha|u|^{2^*(\alpha)}dx
-\frac{\lambda}{p}\int_{\mathbb{R}^N}|x|^\beta|u|^pdx\nonumber\\
&=&\left(\frac{1}{2}-\frac{1}{2}\right)\|u\|^2+\left(\frac{1}{2}-\frac{1}{2^*(\alpha)}
\right)\int_{{\mathbb{R}^N}}|x|^\alpha|u|^{2^*(\alpha)}dx
+\left(\frac{\lambda}{2}-\frac{\lambda}{p}\right)\int_{{\mathbb{R}^N}}|x|^\beta|u|^pdx\nonumber\\
&\geqslant &0.  \label{3.19}
\end{eqnarray}
Let  $v_n:=u_n-u$. By the Brezis-Lieb lemma \cite{1983B-Lieb}, we have
\begin{eqnarray}
\int_{\mathbb{R}^N}|x|^\alpha|u_n|^{2^*(\alpha)}dx-\int_{\mathbb{R}^N}|x|^\alpha|v_n|^{2^*(\alpha)}dx\rightarrow\int_{\mathbb{R}^N} |x|^\alpha|u|^{2^*(\alpha)}dx \ \ \mbox{as} \ n\rightarrow \infty,\label{3.20}
\end{eqnarray}
$$\int_{\mathbb{R}^N}|x|^\beta|u_n|^pdx-\int_{\mathbb{R}^N}|x|^\beta|v_n|^pdx
\rightarrow\int_{\mathbb{R}^N} |x|^\beta|u|^p dx  \ \ \mbox{as} \ n\rightarrow\infty.$$
By $\langle \Phi'(u_n), u_n\rangle \rightarrow 0$  and \eqref{3.18}, we have
\begin{eqnarray*}
&& \|v_n\|^2-\lambda\int_{\mathbb{R}^N}|x|^\beta|v_n|^pdx-\int_{\mathbb{R}^N}|x|^\alpha|v_n|^{2^*(\alpha)}dx\\ &\rightarrow&
-\|u\|^2+\lambda\int_{\mathbb{R}^N}|x|^\beta|u|^pdx+\int_{\mathbb{R}^N}|x|^\alpha|u|^{2^*(\alpha)}\\
&&=-\langle \Phi'(u), u\rangle = 0 \ \hbox{(using \eqref{3.18})}.
\end{eqnarray*}
It follows from $\int_{\mathbb{R}^N}|x|^\beta|v_n|^pdx\rightarrow 0$ under the assumption \eqref{3.8} that
$$\|v_n\|^2-\int_{\mathbb{R}^N}|x|^\alpha|v_n|^{2^*(\alpha)}dx\rightarrow0.$$
Thus we may  assume that
\begin{eqnarray}\label{3.21}
\|v_n\|^2\rightarrow \zeta, \ \   \int_{\mathbb{R}^N}|x|^\alpha|v_n|^{2^*(\alpha)}dx\rightarrow \zeta, \ \ \ n\to\infty.
\end{eqnarray}
By H\'enon-Sobolev inequality \eqref{2.6} we have
$$\|v_n\|^2\geqslant\int_{\mathbb{R}^N}|\nabla v_n|^2dx\geqslant S_\alpha \left(\int_{\mathbb{R}^N}|x|^\alpha|v_n|^{2^*(\alpha)}dx\right)^{\frac{2}{2^*(\alpha)}}.$$
This implies that $\zeta \geqslant S_\alpha\zeta^{\frac{2}{2^*(\alpha)}}$, and so either $\zeta=0$ or $\zeta\geqslant
S_\alpha^{\frac{2+\alpha}{N+\alpha}}.$

 Assume that $\zeta\geqslant S_\alpha^{\frac{N+\alpha}{2+\alpha}}$. By $\Phi(u_n)\rightarrow c$ as $n\rightarrow\infty$, \eqref{3.17} and \eqref{3.20}, we have that
\begin{equation}\label{3.22}
\Phi(u)+\frac{1}{2}\|v_n\|^2-\frac{1}{2^*(\alpha)}
\int_{\mathbb{R}^N}|x|^\alpha|v_n|^{2^*(\alpha)}dx\rightarrow c, \ \ \ \ n\to\infty.
\end{equation}
It follows from \eqref{3.19}, \eqref{3.22} that
\begin{eqnarray*}
 c \geqslant \left(\frac{1}{2}-\frac{1}{2^*(\alpha)}\right)\zeta \geqslant \frac{N+2}{2(N+\alpha)} S_\alpha^{\frac{N+\alpha}{2+\alpha}},
\end{eqnarray*}
which is a contradiction to the assumption that $c<\frac{N+2}{2(N+\alpha)} S_\alpha^{\frac{N+\alpha}{2+\alpha}}$. It must be $\zeta=0$ and then the proof is complete.  \hfill $\Box$

\begin{lem}\label{lem2} Under the assumptions of Theorem {\rm\ref{thm5}}, we have
\begin{eqnarray}\label{3.23}
 m<\frac{2+\alpha}{2(N+\alpha)} S_\alpha^{\frac{N+\alpha}{2+\alpha}}.
\end{eqnarray}
\end{lem}
{\bf Proof.}  Due to \eqref{3.14}, we will get \eqref{3.23} by find a nonzero function $v \in H_r^1(\mathbb{R}^N)$ such that
 \begin{eqnarray}\label{3.24}
 \max_{t \geqslant 0} \Phi(t v)< \frac{2+\alpha}{2(N+\alpha)} S_\alpha^{\frac{N+\alpha}{2+\alpha}}.
 \end{eqnarray}
We first treat the case with the assumption \eqref{3.5} which is contained in $2_*(\beta)<p<2^*(\beta)$ with $\beta>-2$.  In this case we choose $v$ to be the function $u_{\epsilon, \alpha}$ defined by \eqref{2.7}. By the estimations \eqref{2.8}--\eqref{2.10}, we have
\begin{eqnarray}
\lim_{\epsilon\rightarrow0^+}\int_{\mathbb{R}^N}\left|\nabla u_{\epsilon,\alpha}\right|^2dx=S_\alpha^{\frac{N+\alpha}{2+\alpha}},\label{3.25}\\
\lim_{\epsilon\rightarrow0^+}\int_{\mathbb{R}^N}|x|^\alpha |u_{\epsilon,\alpha}|^{2^*(\alpha)}dx=S_\alpha^{\frac{N+\alpha}{2+\alpha}},\label{3.26}\\
\lim_{\epsilon\rightarrow0^+}\int_{\mathbb{R}^N}|u_{\epsilon,\alpha}|^2dx=0.\label{3.27}
\end{eqnarray}
 As $2_*(\beta)<p<2^*(\beta)$,  it follows from \eqref{2.11} that
\begin{eqnarray}
 \lim_{\epsilon\rightarrow0}\int_{\mathbb{R}^N}|x|^\beta|u_{\epsilon,\alpha}|^pdx=0.\label{3.28}
\end{eqnarray}
Since $p>2$,  there is a unique $t_\epsilon>0$ such that
\begin{eqnarray} \label{3.29}
 \sup_{t\geqslant0} \Phi(t u_{\epsilon,\alpha})=\Phi(t_\epsilon u_{\epsilon,\alpha}).
 \end{eqnarray}
It follows that
\begin{eqnarray}
 \|u_{\epsilon,\alpha}\|^2 -t_\epsilon^{2^*(\alpha)-2}\int_{\mathbb{R}^N}|x|^\alpha|u_{\epsilon,\alpha}|^{2^*(\alpha)}dx
-\lambda t_\epsilon^{p-2}\int_{\mathbb{R}^N}|x|^\beta |u_{\epsilon,\alpha}|^pdx=0.\label{3.30}
\end{eqnarray}
Therefore for any $\lambda>0$
\begin{eqnarray*}
0<t_\epsilon \leqslant \left(\frac{\|u_{\epsilon,\alpha}\|^2 }{\int_{\mathbb{R}^N}|x|^\alpha|u_{\epsilon,\alpha}|^{2^*(\alpha)}dx}\right)
^{\frac{1}{2^*(\alpha)-2}}.
\end{eqnarray*}
By \eqref{3.25}--\eqref{3.27} we obtain
$$\lim_{\epsilon\rightarrow0^+}t_\epsilon \leqslant 1, $$
and combining with \eqref{3.30} and \eqref{3.28}, we  have
\begin{eqnarray}
\lim_{\epsilon\rightarrow0^+}t_\epsilon=1.\label{3.31}
\end{eqnarray}
Now \begin{equation} \label{3.32}
\begin{array}{ll}
& \displaystyle \sup_{t\geqslant 0} \Phi (t u_{\epsilon,\alpha}) = \Phi (t_\epsilon u_{\epsilon,\alpha}) \\[3mm]
  \leqslant & \displaystyle \max_{t\geqslant0}\left\{\frac{t^2}{2}\int_{{\mathbb{R}^N}}|\nabla u_{\epsilon,\alpha}|^2-\frac{t^{2^*(\alpha)}}{2^*(\alpha)}
 \int_{{\mathbb{R}^N}}|x|^\alpha|u_{\epsilon,\alpha}|^{2^*(\alpha)}dx\right\}   \\[3mm] & \ \ +\displaystyle \frac{t_\epsilon^2}{2}\int_{{\mathbb{R}^N}}| u_{\epsilon,\alpha}|^2
 -\lambda\frac{ t_\epsilon^p}{p}\int_{{\mathbb{R}^N}}|x|^\beta|u_{\epsilon,\alpha}|^pdx \\[3mm]
 = &  \displaystyle \frac{2+\alpha}{2(N+\alpha)}S_\alpha^{\frac{N+\alpha}{2+\alpha}}+O\left(\epsilon^{N-2}\right)
 +\frac{t_\epsilon^2}{2}\int_{{\mathbb{R}^N}}| u_{\epsilon,\alpha}|^2-\lambda \frac{ t_\epsilon^p}{p}\int_{{\mathbb{R}^N}}|x|^\beta|u_{\epsilon,\alpha}|^pdx.  \end{array}
 \end{equation}
 We use \eqref{2.10} and \eqref{2.11} to prove
 \begin{equation} \label{3.33}
 O\left(\epsilon^{N-2}\right)
 +\frac{t_\epsilon^2}{2}\int_{{\mathbb{R}^N}}| u_{\epsilon,\alpha}|^2-\lambda \frac{ t_\epsilon^p}{p}\int_{{\mathbb{R}^N}}|x|^\beta|u_{\epsilon,\alpha}|^pdx<0 \ \ \hbox{for} \ \epsilon>0 \ \hbox{small}.\end{equation}
   For the case of \eqref{3.5} with $N=3$, \ $3+\beta-\frac p2<1<\frac p2$,
 $$  \int_{{\mathbb{R}^N}}| u_{\epsilon,\alpha}|^2 dx = O(\epsilon), \ \  \int_{{\mathbb{R}^N}}|x|^\beta|u_{\epsilon,\alpha}|^pdx = \left \{\begin{array}{ll}
 O(\epsilon^{\frac{3+\beta}{2}})+C\epsilon^{\frac{3+\beta}{2}}|\ln \epsilon|, &   p=3+\beta;\\
 O(\epsilon^{\frac{p}{2}})+C\epsilon^{3+\beta-\frac{p}{2}}, &   p>3+\beta.
\end{array}
\right.  $$
For the case of \eqref{3.5} with $N=3$, \ $3+\beta-\frac p2<1<\frac p2$,
\begin{eqnarray*}
\int_{{\mathbb{R}^N}}| u_{\epsilon,\alpha}|^2 dx = O(\epsilon), \ \  \int_{{\mathbb{R}^N}}|x|^\beta|u_{\epsilon,\alpha}|^pdx =
 O(\epsilon^{\frac{p}{2}})+C\epsilon^{3+\beta-\frac{p}{2}},   p>2>3+\beta.
\end{eqnarray*}
For the case of \eqref{3.5} with $N=4$,  \ $4+\beta-p<2<p$,
  $$  \int_{{\mathbb{R}^N}}| u_{\epsilon,\alpha}|^2 dx =
 O(\epsilon^2)+C\epsilon^2|\ln \epsilon|,$$ $$
\int_{\mathbb{R}^N}|x|^\beta|u_{\epsilon,\alpha}|^pdx=\left \{\begin{array}{ll}
 O(\epsilon^{\frac{4+\beta}{2}})+C\epsilon^{\frac{4+\beta}{2}}|\ln \epsilon|, &   p=\frac{4+\beta}{2};\\
 O(\epsilon^{p})+C\epsilon^{4+\beta-p},&   p>\frac{4+\beta}{2}.
\end{array}\right.$$
For the case of \eqref{3.5} with $N\geqslant 5$, \ $N+\beta-\frac{p(N-2)}{2}<2<\frac{p(N-2)}{2}$,
   $$  \int_{{\mathbb{R}^N}}| u_{\epsilon,\alpha}|^2 dx =
 O(\epsilon^{N-2})+C\epsilon^2, $$
$$
\int_{\mathbb{R}^N}|x|^\beta|u_{\epsilon,\alpha}|^pdx=\left \{\begin{array}{ll}
 O(\epsilon^{\frac{N+\beta}{2}})+C\epsilon^{\frac{N+\beta}{2}}|\ln \epsilon|, &   p=\frac{N+\beta}{N-2};\\
 O(\epsilon^{\frac{p(N-2)}{2}})+C\epsilon^{N+\beta-\frac{p(N-2)}{2}}, &   p>\frac{N+\beta}{N-2},
\end{array}
\right.$$
 It follows that \eqref{3.33} holds for $\epsilon>0$ small enough.
  Therefore by \eqref{3.32}, and (\ref{3.33}), under the assumption \eqref{3.5} we have that for $\epsilon>0$ small enough
\begin{eqnarray}
\sup_{t\geqslant0}\Phi(t u_{\epsilon,\alpha}) = \Phi(t_\epsilon u_{\epsilon,\alpha})<\frac{2+\alpha}{2(N+\alpha)}S_\alpha^{\frac{N+\alpha}{2+\alpha}} \label{3.34}
\end{eqnarray} for any $\lambda>0$.

Secondly, we consider $p=2_*(\beta)=2$ and $-2<\beta<0, 0<\lambda<\lambda_{1\beta}$. Applying the inequality \eqref{3.11}, it is easy to see that \eqref{3.31} and \eqref{3.32} hold with $p=2$. We now prove \eqref{3.33} is still true under the assumptions of (iii) in Theorem \ref{thm5}.  As $N=3, -2<\beta\leq -1$,
 $$  \int_{{\mathbb{R}^N}}| u_{\epsilon,\alpha}|^2 dx = O(\epsilon), \ \  \int_{{\mathbb{R}^N}}|x|^\beta|u_{\epsilon,\alpha}|^2dx = \left \{\begin{array}{ll}
 O(\epsilon)+C\epsilon|\ln \epsilon|, &   \beta=-1;\\
 O(\epsilon)+C\epsilon^{2+\beta}, &   -2<\beta<-1.
\end{array}
\right.  $$
As $N=4, -2<\beta<0$, then
\begin{eqnarray*}
 \int_{{\mathbb{R}^N}}| u_{\epsilon,\alpha}|^2 dx=O(\epsilon^2)+C\epsilon^2|\ln \epsilon|,\
 \int_{\mathbb{R}^N}|x|^\beta|u_{\epsilon,\alpha}|^2dx=
 O(\epsilon^2)+C\epsilon^{2+\beta}.
\end{eqnarray*}
As $N\geq5, -2<\beta<0$, then
\begin{eqnarray*}
\int_{{\mathbb{R}^N}}| u_{\epsilon,\alpha}|^2 dx=O(\epsilon^{N-2})+C\epsilon^2,\
\int_{\mathbb{R}^N}|x|^\beta|u_{\epsilon,\alpha}|^2dx=
 O(\epsilon^{N-2})+C\epsilon^{2+\beta}.
\end{eqnarray*}
Thus under (iii) in Theorem \ref{thm5}, we know \eqref{3.34} is true for $\epsilon>0$ small enough.

We next consider the case that $2_*(\beta)<p<2^*(\beta)$ with $\beta>-2$. We choose a function
  $\phi\in C^\infty_{0,r}({\mathbb{R}^N})$ such that $\phi(|x|)\geqslant 0$ and $\phi(0)=1$. Set $v_0(x)=\phi(|x|)|x|^{-k}$ with $k\in (0, \frac12)$ so that $v_0\in H_r^1(\mathbb{R}^N)$ and $\|v_0\|_{L^{2^*(\alpha)}({\mathbb{R}^N};|x|^\alpha)}>0$.

  Since $p>2$, it is easy to see that $\sup_{t \geqslant 0}\Phi(t v_0)$ is achieved at a unique $t_\lambda>0$  for each $\lambda>0$. Then
\begin{eqnarray} \label{3.35}
\sup_{t\geqslant0}\Phi(t v_0)&=&\frac{t_{\lambda}^2}{2} \|v_0\|^2  -\frac{t_\lambda^{2^*(\alpha)}}{2^*(\alpha)}\int_{\mathbb{R}^N}|x|^\alpha|v_0|^{2^*(\alpha)}dx
-\frac{\lambda t_{\lambda}^p}{p}\int_{\mathbb{R}^N}|x|^\beta|v_0|^pdx,
\end{eqnarray}
\begin{equation}
\|v_0\|^2
-t_{\lambda}^{2^*(\alpha)-2}\int_{\mathbb{R}^N}
|x|^\alpha|v_0 |^{2^*(\alpha)}dx-\lambda t_{\lambda}^{p-2}
\int_{\mathbb{R}^N}|x|^\beta|v_0 |^pdx=0.\label{3.36}
\end{equation}
By \eqref{3.36} we get
\begin{eqnarray}
t_\lambda \leqslant  \left(\frac{\|v_0\|^2}{\lambda\int_{\mathbb{R}^N}|x|^\beta|v_0|^pdx}\right)^{\frac{1}{p-2}}.\label{3.37}
\end{eqnarray}
It follows from \eqref{3.36} and \eqref{3.37} that
\begin{eqnarray}
\lim_{\lambda\rightarrow\infty}t_\lambda =0.\label{3.38}
\end{eqnarray}
By \eqref{3.35} and \eqref{3.38} we have
\begin{eqnarray*}
\lim_{\lambda\rightarrow \infty}   \Phi(t_\lambda v_0)\leqslant 0.
\end{eqnarray*}
Thus there is $\lambda^*>0$ such that for all $\lambda>\lambda_*$,
\begin{eqnarray}
\sup_{t\geqslant 0} \Phi (tv_0)= \Phi (t_\lambda v_0) <\frac{2+\alpha}{2(N+\alpha)}S_\alpha^{\frac{N+\alpha}{2+\alpha}}.\label{3.39}
\end{eqnarray}
Therefore, in both cases we get the existence of  $v^* \in H_r^1(\mathbb{R}^N)\backslash \{0\}$ satisfying \eqref{3.24}. The proof is complete. \hfill$\Box$

\noindent {\bf Proof of Theorem \ref{thm5}} \ By Proposition \ref{pro1} and Ekeland's variational principle, there exists a minimizing sequence $\{u_n\}\subset \mathcal{N}$ such that
 \begin{eqnarray}\label{3.40}
 \Phi(u_n)\rightarrow m, \ \ \ \ (\Phi|_{\mathcal{N}})'(u_n)\rightarrow0 \ \ \ n\to\infty.
 \end{eqnarray}
 Let $\lambda_n$ be the Lagrange multiplier satisfying
\begin{eqnarray}
(\Phi|_{\mathcal{N}})'(u_n)=\Phi'(u_n)-\lambda_n \Psi'(u_n).\label{3.41}
\end{eqnarray}
Similar with \eqref{3.15} and \eqref{3.16}, we get that the sequence $\{u_n\}$ is bounded, which implies
that $\Psi'(u_n)$ is bounded. Combining with \eqref{3.41}, one has
 $$\|(\Phi|_{\mathcal{N}})'(u_n)\| \rightarrow0 \ \ \ n\to\infty.$$
Hence
\begin{eqnarray}
o(1)=\langle \Phi'(u_n),u_n\rangle-\lambda_n\langle \Psi'(u_n), u_n\rangle \ \ \ n\to\infty.\label{3.42}
\end{eqnarray}
Since $u_n\in \mathcal{N}$,  $\langle \Phi'(u_n), u_n\rangle=0$. The arguments of the proof in Proposition \ref{pro1}  implies $|\langle \Psi'(u_n), u_n\rangle|>0$. It follows from  \eqref{3.42} that
\begin{eqnarray}
\lambda_n\rightarrow0 \ \  \ \ n\to\infty.\label{3.43}
\end{eqnarray}
Since $\Psi'(u_n)$ is bounded, by \eqref{3.41} and \eqref{3.43}, one has $\Phi'(u_n)\rightarrow 0$ as $n\rightarrow\infty$. Therefore   $\{u_n\}$ is a ${{\rm(PS)}_m}$ sequence of $\Phi$ in $H_r^1(\mathbb{R}^N)$. The boundedness of $\{u_n\}$ in $H_r^1(\mathbb{R}^N)$, Lemma \ref{lem1} and  Lemma \ref{lem2} imply that there exists $u_0\in H_r^1(\mathbb{R}^N)$ such that $$\Phi(u_0)=m, \ \ \ \ \Phi'(u_0)=0.$$ It is clear that $u_0$ is nontrivial.

 It is easy to see that $|u_0| \in \mathcal{N}$ and $\Phi(|u_0|)=m$. By the Lagrange multiplier theorem, there exists a $\lambda\in \mathbb{R}$ such that
 $$\Phi'(|u_0|)=\lambda \Psi'(|u_0|).$$
Thus
\begin{eqnarray*}
0=\langle \Phi'(|u_0|),|u_0|\rangle=\lambda\langle \Psi'(|u_0|),|u_0|\rangle.
\end{eqnarray*}
It follows from the regularity of $\mathcal{N}$ that
$$\lambda=0.$$
Hence
$$\Phi'(|u_0|)=0$$
and $|u_0|$ is a critical point of $\Phi$. By the strong maximum principle, we have $u_0>0$. \hfill$\Box$


\section{Problems with double critical exponents}
In this section we consider the following equation with double critical exponents
\begin{eqnarray}
\left \{\begin{array}{ll}
\displaystyle -\Delta u+u=|x|^{\alpha_1}|u|^{2^*(\alpha_1)-2}u+\mu|x|^{\alpha_2} |u|^{2^*(\alpha_2)-2}u+\lambda |x|^{\beta} |u|^{p-2}u & \hbox{in} \ {\mathbb{R}^N},\\
\displaystyle u\in H_r^1({\mathbb{R}^N}),
\end{array}\label{4.1}
\right.
\end{eqnarray}
where $N \geqslant 3$, $\alpha_1>\alpha_2>-2$, $\beta>-2,$ $\mu\in \mathbb{R}$ and  $\lambda>0$ are parameters.  The corresponding energy functional of \eqref{4.1} reads as
\begin{equation}\label{4.2} \begin{array}{ll}
\Phi_{\lambda, \mu}(u)= & \displaystyle \frac{1}{2}\int_{\mathbb{R}^N}(|\nabla u|^2+|u|^2) dx
-\frac{1}{2^*(\alpha_1)}\int_{\mathbb{R}^N}|x|^{\alpha_1} |u|^{2^*(\alpha_1)}dx  \\[3mm] & \displaystyle -\frac{\mu}{2^*(\alpha_2)}\int_{\mathbb{R}^N}|x|^{\alpha_2}|u|^{2^*(\alpha_2)}dx
-\frac{\lambda}{p}\int_{\mathbb{R}^N}|x|^{\beta}|u|^pdx. \end{array}
\end{equation}
 It is obvious that $\Phi_{\lambda, \mu} \in C^2(H_r^1(\mathbb{R}^N), \mathbb{R})$.
 The Nehari manifold of $\Phi_{\lambda, \mu}$ is defined as
 \begin{eqnarray} \label{4.3}
 \mathcal{N}_{\lambda, \mu}:=\{u\in H^1_r(\mathbb{R}^N) \setminus\{0\}: \ \ \Psi_{\lambda, \mu}(u):= \langle \Phi_{\lambda, \mu}'(u),u\rangle=0\}.
 \end{eqnarray}
 Define
\begin{eqnarray} \label{4.4}
  \displaystyle m_{\lambda, \mu}:=\inf_{u \in \mathcal{N}_{\lambda, \mu}} \Phi_{\lambda, \mu} (u).
\end{eqnarray}
We will prove that in some suitable situations the minimum $m_{\lambda, \mu}$ can be achieved so that \eqref{4.1} has nonnegative  ground state solutions in $H_r^1(\mathbb{R}^N)$. We first give two radial inequalities for functions in $D_r^{1,2}(\mathbb{R}^N)$.
\begin{lem}\label{lem4.1}
Assume $N \geqslant 3$, $\sigma>\varsigma>-2$. For any $u\in D^{1,2}_r(\mathbb{R}^N)$, it holds that
\begin{eqnarray*}
\int_{\mathbb{R}^N}|x|^\sigma|u|^{2^*(\sigma)} dx \leqslant \tilde{C} \|\nabla u\|_{L^2(\mathbb{R}^N)}^{2^*(\sigma)-2^*(\varsigma)}
\int_{\mathbb{R}^N}|x|^\varsigma|u|^{2^*(\varsigma)}dx,
\end{eqnarray*}
where $ {\widetilde{C}}=\left[(N-2)\omega_N\right]^{\frac{\varsigma-\sigma}{N-2}}$.
\end{lem}
{\bf Proof} \ For any $u\in D_r^{1,2}(\mathbb{R}^N)$, we have by \cite[Lemma 1]{2007SWW-1,2007SWW-2} that
\begin{eqnarray}
|u(|x|)|\leqslant \hat{C}|x|^{-\frac{N-2}{2}}\|\nabla u\|_{L^2(\mathbb{R}^N)},\label{19}
\end{eqnarray}
where $\hat{C}=\omega_N^{-\frac{1}{2}}\left(\frac{1}{N-2}\right)^{\frac{1}{2}}$. Since $\sigma>\varsigma$, we have
\begin{eqnarray*}
\int_{\mathbb{R}^N}|x|^\sigma|u|^{2^*(\sigma)} dx &=&\int_{\mathbb{R}^N}|x|^\varsigma|x|^{\sigma-\varsigma}
|u|^{2^*(\varsigma)}|u|^{2^*(\sigma)-2^*(\varsigma)} dx\\
& \leqslant &\tilde{C} \|\nabla u\|_{L^2(\mathbb{R}^N)}^{2^*(\sigma)-2^*(\varsigma)}
\int_{\mathbb{R}^N}|x|^\varsigma|u|^{2^*(\varsigma)}dx,
\end{eqnarray*}
where $\tilde{C}=\hat{C}^{2^*(\sigma)-2^*(\varsigma)}=\left[(N-2)\omega_N\right]^{\frac{\varsigma-\sigma}{N-2}}$. $\hfill\Box$
\begin{lem}\label{lem4.2}
Assume $N \geqslant 3, \varsigma>\sigma>-2$. For any $u\in D_r^{1,2}(\mathbb{R}^N)$, it holds that \begin{eqnarray*}
\|u\|_{L^{2^*(\sigma)}(\mathbb{R}^N;|x|^\sigma)}\leqslant S_\theta^{-\frac{1-\tau}{2}}\|u\|_{L^{2^*(\varsigma)}(\mathbb{R}^N;|x|^\varsigma)}^\tau\|\nabla u\|_{L^2(\mathbb{R}^N)}^{1-\tau},
\end{eqnarray*}
where $\theta= \frac{ 2^*(\varsigma)\sigma -m \varsigma}{2^*(\varsigma)-m}$, $\tau=\frac{m}{2^*(\sigma)}\in\left(0,\frac{(2+\sigma)(N+\varsigma)}
{(2+\varsigma)(N+\sigma)}\right]$,  $0<m\leqslant\frac{2+\sigma}{2+\varsigma}2^*(\varsigma). $
\end{lem}
{\bf Proof} \ For $u\in D_r^{1,2}(\mathbb{R}^N)$, we choose $0<m\leqslant \frac{2+\sigma}{2+\varsigma}2^*(\varsigma)$. By H\"{o}lder inequality, we get
\begin{eqnarray*}
&&\int_{\mathbb{R}^N}|x|^\sigma|u|^{2^*(\sigma)}dx\\ &\leqslant&\left(\int_{\mathbb{R}^N}|x|^\varsigma|u|^{2^*(\varsigma)} dx \right)^{\frac{m}{2^*(\varsigma)}}
\left(\int_{\mathbb{R}^N}|x|^{\left(\sigma-\frac{m\varsigma}{2^*(\varsigma)}\right)\frac{2^*(\varsigma)}
{2^*(\varsigma)-m}}|u|^{2^*(\varsigma)\frac{2^*(\sigma)-m}{2^*(\varsigma)-m}}
 dx \right)^{\frac{2^*(\varsigma)-m}{2^*(\varsigma)}}.
\end{eqnarray*}
Set $\theta=\left(\sigma-\frac{m\varsigma}{2^*(\varsigma)} \right)\frac{2^*(\varsigma)}
{2^*(\varsigma)-m}.$ Then \begin{eqnarray*}
2^*(\theta)=\frac{2^*(\sigma)-m}{2^*(\varsigma)-m}2^*(\varsigma).
\end{eqnarray*}
It follows that
\begin{eqnarray*}
\|u\|_{L^{2^*(\sigma)}(\mathbb{R}^N;|x|^\sigma)} \leqslant  S_\theta^{-\frac{1-\tau}{2}}\|u\|_{L^{2^*(\varsigma)}(\mathbb{R}^N;|x|^\varsigma)}^\tau\|\nabla u\|_{L^2(\mathbb{R}^N)}^{1-\tau},
\end{eqnarray*}
where $\tau=\frac{m}{2^*(\sigma)}\in\left(0,\frac{(2+\sigma)(N+\varsigma)}{(2+\varsigma)(N+\sigma)}\right], 0< m \leqslant \frac{2+\sigma}{2+\varsigma}2^*(\varsigma)$ and $S_\theta$ is the embedding constant from $D_r^{1,2}(\mathbb{R}^N)$ into $L^{2^*(\theta)}(\mathbb{R}^N; |x|^\theta)$. \hfill$\Box$

 \subsection{The Case $\mu>0$}
 In this subsection,  we establish the existence of ground state solution for \eqref{4.1} for the case $\mu>0$.  We will prove the following theorem.
 \begin{thm}\label{thm6}
 Assume $2_*(\beta)<p<2^*(\beta)$ with $\beta>-2$. For any $\mu>0$ being fixed, there exists $\lambda^*>0$
 such that the equation \eqref{4.1} has a positive ground state solution for $\lambda>\lambda^*$.
 \end{thm}

 We may assume $\mu=1$ and set $\Phi_{\lambda,\mu}=\Phi_\lambda, \mathcal{N}_{\lambda,\mu}=\mathcal{N}_\lambda, m_{\lambda,\mu}=m_\lambda$. First we have the same properties on the corresponding
 Nehari manifold.

  \begin{pro}\label{pro4.4}  Assume $\beta>-2$ and $2_*(\beta)<p<2^*(\beta)$.
  \begin{enumerate}
 \item[{\rm (i) }] For $u\in H_r^1(\mathbb{R}^N)\backslash\{0\}$, there exists a unique $t_u>0$ such that  $t_u u\in\mathcal{N}_\lambda$ and $\Phi_\lambda (t_u u)=\max_{t\geqslant  0}\Phi_\lambda (tu)$, moreover, the manifold $\mathcal{N}_\lambda$ is nonempty;
     \item[{\rm (ii) }] The manifold $\mathcal{N}_\lambda$ is $C^1$ regular, that is $\langle \Psi_\lambda'(u), u\rangle=\langle \Phi_\lambda''(u)u, u\rangle\ne0$ for any $u\in \mathcal{N}_\lambda$;
     \item[{\rm (iii) }] $\mathcal{N}_\lambda$ is closed and bounded away from $0$, and $m_\lambda>0$.
 \end{enumerate}
  \end{pro}

Following the arguments of \cite[Chapter 4]{1996Willem}, we can prove that
 \begin{eqnarray}  \label{4.6}
        m_\lambda = c_\lambda:=\inf_{u\in H_r^1(\mathbb{R}^N) \setminus\{0\}}
  \max_{t\geqslant 0} \Phi_\lambda (t u).
 \end{eqnarray}

Next we verify the  local (PS) condition for $\Phi_\lambda$. We have

\begin{lem}\label{lem4.5}  Any a $(PS)_c$ sequence $\{u_n\}\subset H^1_r(\mathbb{R}^N)$ of $\Phi_\lambda$ with $$c\leqslant c^*:=\frac{2+\alpha_2}{2(N+\alpha_2)} {\widetilde{M}}$$ contains a convergent subsequence, where $ {\widetilde{M}}:=M(\alpha_1,\alpha_2,N,S_{\alpha_2})$ is the unique positive solution of the equation
 \begin{eqnarray} \label{403}
 \widetilde{C}t^{\frac{2^*(\alpha_1)-2}{2}}+t^{\frac{2^*(\alpha_2)-2}{2}}-S_{\alpha_2}
 ^{\frac{2^*(\alpha_2)}{2}}=0
 \end{eqnarray}
 and $\widetilde{C}$ is given in Lemma {\rm \ref{lem4.1}.}
 \end{lem}
{\bf Proof.} \ Let $\{u_n\}\subset H^1_r({\mathbb{R}^N})$ be such that
 \begin{eqnarray} \label{4034} \Phi_\lambda (u_n) \to c, \ \ \ \ \ \Phi_\lambda'(u_n) \to 0 \ \ \ \ n\to\infty.\end{eqnarray}
Take $\tilde{p}:=\min\{2^*(\alpha_2),p\}.$  For $n$ large enough, we have
\begin{eqnarray*}
&&c+1+o(1) \|u_n\| \\
&\geqslant& \Phi_\lambda(u_n)-\frac{1}{\tilde{p}}\langle \Phi_\lambda'(u_n),u_n\rangle\\
&=&\left(\frac{1}{2}-\frac{1}{\tilde{p}}\right)\| u_n\|^2
+\sum_{i=1}^2\left(\frac{1}{\tilde{p}}-\frac{1}{2^*(\alpha_i)}\right)
\int_{{\mathbb{R}^N}}|x|^{\alpha_i}|u_n|^{2^*(\alpha_i)}dx\\ && \ \ \ \ +\left(\frac{1}{\tilde{p}}
-\frac{1}{p}\right)\int_{{\mathbb{R}^N}}|x|^\alpha|u_n|^pdx\\
&\geqslant& \left(\frac{1}{2}-\frac{1}{\tilde{p}}\right)\|u_n\|^2.
\end{eqnarray*}
Since $\tilde{p}>2$, it follows that  $\{u_n\}$ is bounded in $H^1_r(\mathbb{R}^N)$. Up to a subsequence, we may  assume that there is $u\in H^1_r(\mathbb{R}^N)$ such that
\begin{eqnarray}
  \left\{\begin{array}{ll}
&u_n\rightharpoonup u \ \ \mbox{in} \ H^1_r(\mathbb{R}^N);\\
&u_n\rightharpoonup u \ \ \mbox{in} \ L^{2^*(\gamma)}(\mathbb{R}^N;|x|^\gamma), \ \gamma>-2;\\
&u_n\rightarrow u \ \  \mbox{in} \ L^p(\mathbb{R}^N;|x|^\beta),  \ \   \ 2_*(\beta)<p<2^*(\beta),  \ \beta>-2;\\
&u_n(x)\rightarrow u(x) \ \ \mbox{a.e.\ on} \  \mathbb{R}^N.
  \end{array}
\right.\label{4009}
\end{eqnarray}
  It follows from \eqref{4009} that for any $w\in H^1_r({\mathbb{R}^N})$
$$\int_{\mathbb{R}^N} |x|^{\alpha_i}|u_n|^{2^*(\alpha_i)-2}u_nwdx\rightarrow \int_{\mathbb{R}^N}|x|^{\alpha_i} |u|^{2^*(\alpha_i)-2}uwdx, \ i=1, 2.$$
 By \eqref{4034} we have that $u$ solves weakly the equation
 $$-\Delta u+u=|x|^{\alpha_1}|u|^{2^*(\alpha_1)-2}u+|x|^{\alpha_2}|u|^{2^*(\alpha_2)-2}u+\lambda |x|^\beta|u|^{p-2}u.$$
Therefore $\langle \Phi_\lambda'(u),u\rangle=0$ and
\begin{eqnarray*}
\Phi_\lambda (u) =  \Phi_\lambda(u)-\frac{1}{\tilde{p}}\langle \Phi_\lambda'(u),u\rangle
\geqslant  \left(\frac{1}{2}-\frac{1}{\tilde{p}}\right)\|u \|^2 \geqslant 0.
\end{eqnarray*}
 Set $v_n:=u_n-u$. By \eqref{4009} and  the Brezis-Lieb lemma(\cite{1983B-Lieb}), we have
$$\int_{\mathbb{R}^N}|x|^{\alpha_i}|u_n|^{2^*(\alpha_i)} dx
-\int_{{\mathbb{R}^N}}|x|^{\alpha_i}|v_n|^{2^*(\alpha_i)} dx
\rightarrow\int_{\mathbb{R}^N} |x|^{\alpha_i}|u|^{2^*(\alpha_i)} dx,\ \mbox{as}\ n\rightarrow\infty,i=1,2.$$
It follows from \eqref{4034} that $\langle \Phi_\lambda'(u_n), u_n\rangle \rightarrow 0$ as $n\to\infty$. Thus
\begin{eqnarray*}
&&\|v_n\|^2
-\sum_{i=1}^2\int_{\mathbb{R}^N}|x|^{\alpha_i}
|v_n|^{2^*(\alpha_i)}dx-\lambda\int_{\mathbb{R}^N}|x|^\beta|v_n|^pdx\\
&\rightarrow& - \|u\|^2
+\sum_{i=1}^2\int_{\mathbb{R}^N}|x|^{\alpha_i}|u|^{2^*(\alpha_i)}dx
+\lambda\int_{\mathbb{R}^N}|x|^\beta|u|^pdx = -\langle \Phi_\lambda'(u),u\rangle = 0,
\end{eqnarray*}
Since
\begin{eqnarray*}
\int_{\mathbb{R}^N}|x|^\beta |v_n|^pdx\rightarrow 0,
\end{eqnarray*}
it follows that
\begin{equation}\label{20}
\|v_n\|^2
-\sum_{i=1}^2\int_{\mathbb{R}^N}|x|^{\alpha_i}|v_n|^{2^*(\alpha_i)}dx\rightarrow0.
\end{equation}
 We assume, up to a subsequence if necessary,  that
$$\displaystyle \lim_{n\rightarrow\infty} \|v_n\| =A_\infty,\ \ \lim_{n\rightarrow\infty} \int_{\mathbb{R}^N}|x|^{\alpha_1}|v_n|^{2^*(\alpha_1)}dx=B_\infty,$$
$$\ \ \lim_{n\rightarrow\infty}
\int_{\mathbb{R}^N}|x|^{\alpha_2}|v_n|^{2^*(\alpha_2)}dx=C_\infty.$$
By (\ref{20}), we have
\begin{eqnarray} \label{406} A_\infty=B_\infty+C_\infty.\end{eqnarray}
We will end the proof by proving $A_\infty=0$.

 Assume that $A_\infty>0$.   By Lemma  \ref{lem4.1} and Lemma \ref{lem4.2}, we have
\begin{eqnarray*}
B_\infty \leqslant  \tilde{C} A_\infty^{\frac{2^*(\alpha_1)-2^*(\alpha_2)}{2}}C_\infty,  \ \ \ \
C_\infty^{\frac{1}{2^*(\alpha_2)}}\leqslant S_\theta^{-\frac{1-\tau}{2}}A_\infty^{\frac{1-\tau}{2}}B_\infty^{\frac{\tau}{2^*(\alpha_1)}},
\end{eqnarray*}
where
 $$\tilde{C}=\left[(N-2)\omega_N\right]^{\frac{\alpha_2-\alpha_1}{N-2}}, \ \ \ \ \theta=\frac{2^*(\alpha_1)\alpha_2- m\alpha_1}{2^*(\alpha_1)-m},$$
$$\tau=\frac{m}{2^*(\alpha_2)}\in \left(0,\frac{(2+\alpha_2)(N+\alpha_1)}
{(2+\alpha_1)(N+\alpha_2)}\right], \ \ \ 0< m \leqslant \frac{2+\alpha_2}{2+\alpha_1}2^*(\alpha_1).$$
It must be $B_\infty>0$ and $C_\infty>0$. By H\'enon-Sobolev inequality, we have
\begin{eqnarray}\label{407}
S_{\alpha_2} C_\infty^{\frac{2}{2^*(\alpha_2)}}\leqslant  A_\infty=B_\infty+C_\infty \leqslant  \tilde{C} A_\infty^{\frac{2^*(\alpha_1)-2^*(\alpha_2)}{2}}C_\infty+C_\infty.
\end{eqnarray}
Thus
\begin{eqnarray} \label{408}
C_\infty\geqslant  {S_{\alpha_2}^{\frac{2^*(\alpha_2)}{2^*(\alpha_2)-2}}}{\left(1+\tilde{C}A_\infty^
{\frac{2^*(\alpha_1)-2^*(\alpha_2)}{2}}\right)^{-\frac{2^*(\alpha_2)}{2^*(\alpha_2)-2}}}.
\end{eqnarray}
By \eqref{407} and \eqref{408} we deduce that
\begin{eqnarray*}
\tilde{C}A_\infty^{\frac{2^*(\alpha_1)-2}{2}}+ A_\infty^{\frac{2^*(\alpha_2)-2}{2}}
-S_{\alpha_2}^{\frac{2^*(\alpha_2)}{2}} \geqslant  0.
\end{eqnarray*}
It follows that $A_\infty \geqslant \widetilde{M}$
where $\widetilde{M}$ is the unique positive solution of the
 equation \eqref{403}.
Since $\Phi_\lambda (u_n)\rightarrow c$ as $n\rightarrow\infty$,
\begin{equation*}
  \Phi_\lambda(u)+\frac{1}{2}\|v_n\|^2
  -\sum_{i=1}^2\frac{1}{2^*(\alpha_i)}\int_{{\mathbb{R}^N}}|x|^{\alpha_i}|v_n|^{2^*(\alpha_i)}dx\rightarrow c.
\end{equation*}
As $\Phi_\lambda(u)\geqslant0$,  we have by \eqref{406} that \begin{eqnarray*}
c  \geqslant   \frac{A_\infty}{2}-\frac{B_\infty}{2^*(\alpha_1)}-\frac{ C_\infty}{2^*(\alpha_2)}
> \frac{2+\alpha_2}{2(N+\alpha_2)} {\widetilde{M}}.
\end{eqnarray*}
It contradicts the choice of $c$.  The proof is complete.  \hfill $\Box$

 \begin{lem}\label{lem4.6}  Assume $2_*(\beta)<p<2^*(\beta)$ with $\beta>-2$.
 There exists $\lambda^*>0$  such that for $\lambda>\lambda^*$,
 \begin{eqnarray*}
 m_\lambda= c_\lambda \leqslant  c^*:=\frac{2+\alpha_2}{2(N+\alpha_2)} \widetilde{M},
 \end{eqnarray*}
 where ${\widetilde{M}}>0$ is the unique positive solution of \eqref{403}.
 \end{lem}
 {\bf Proof} \ We choose a function
  $\phi\in C^\infty_{0,r}({\mathbb{R}^N})$ such that $\phi(|x|)\geqslant 0$ and $\phi(0)=1$. Set $v_0(x)=\phi(|x|)|x|^{-k}$ with $k\in (0, \frac12)$ so that $v_0\in H_r^1(\mathbb{R}^N)$ and $\|v_0\|_{L^{2^*(\alpha_i)}({\mathbb{R}^N};|x|^{\alpha_i})}>0,i=1,2$. It is easily seen  that $\displaystyle\sup_{t\geqslant0}\Phi_\lambda(t v_0)$ is achieved at a unique $t_\lambda>0$ so that
 \begin{eqnarray} \label{4.14}
\sup_{t\geqslant 0}\Phi_\lambda(t v_0 ) =  \frac{t_\lambda^2}{2} \|v_0\|^2 -\sum_{i=1}^2\frac{t_\lambda^{2^*(\alpha_i)}}{2^*(\alpha_i)}\int_{\mathbb{R}^N}
 |x|^{\alpha_i}|v_0|^{2^*(\alpha_i)}dx -\frac{\lambda t_\lambda^p}{p}\int_{\mathbb{R}^N}|x|^\beta|v_0|^pdx,
 \end{eqnarray}
 and
 \begin{eqnarray*}
 \|v_0\|^2
 = \sum_{i=1}^2t_\lambda^{2^*(\alpha_i)-2}\int_{\mathbb{R}^N}
 |x|^{\alpha_i}|v_0|^{2^*(\alpha_i)}dx+\lambda t_{\lambda}^{p-2}\int_{\mathbb{R}^N}|x|^\beta|v_0|^pdx.
 \end{eqnarray*} It follows that
 \begin{eqnarray*}
 0< t_\lambda \leqslant \left(\frac{\|v_0\|^2} {\lambda \int_{\mathbb{R}^N}|x|^\beta|v_0|^pdx}\right)
 ^{\frac{1}{p-2}}.
 \end{eqnarray*}
 Hence
\begin{eqnarray}
\lim_{\lambda\rightarrow\infty}t_\lambda =0.\label{4.15}
\end{eqnarray}
By \eqref{4.14} and \eqref{4.15} we have that
  $$\lim_{\lambda\to \infty} \sup_{t\geqslant 0} \Phi_\lambda(t v_0) \leqslant 0.$$
 Therefore there exists a $\lambda^*>0$ such that for $\lambda>\lambda^*$
\begin{eqnarray*}
m_\lambda =c_\lambda=\inf_{u\in H_r^1(\mathbb{R}^N), u\ne0}\sup_{t\geqslant 0} \Phi_\lambda(t u) \leqslant \sup_{t\geqslant 0} \Phi_\lambda(t v_0) \leqslant \frac{2+\alpha_2}{2(N+\alpha_2)} {\widetilde{M}}.
\end{eqnarray*}
The proof is complete. \hfill$\Box$

\noindent {\bf Proof  of Theorem \ref{thm6}} \ The argument is same as that of Theorem \ref{thm5}. By Proposition \ref{pro4.4} and Ekeland's variational principle, there exists a minimizing sequence $\{u_n\}\subset \mathcal{N}_\lambda$ such that $\Phi_\lambda (u_n)\rightarrow m_\lambda, (\Phi_\lambda|_{\mathcal{N}_\lambda})'(u_n)\rightarrow0$.  The boundedness of $\{u_n\}$ implies that $\Psi_\lambda'(u_n)$ is bounded,   Proposition \ref{pro4.4}, Lemma \ref{lem4.5} and  Lemma \ref{lem4.6} imply that there exists nonnegative $u \in H_r^1(\mathbb{R}^N) $ such that $\Phi_\lambda(u)=m_\lambda$ and $\Phi_\lambda'(u)=0$. By the strong maximum principle, we have $u>0$. \hfill$\Box$

 \subsection{The case $\mu<0$}

 In this subsection,  we establish the existence of ground state solution for \eqref{4.1} for the case $\mu<0$ and $\lambda>0$. In this case the parameter $\mu$ plays some role and the range of the power $p$ is more restrictive.
 We will prove the following theorem.
 \begin{thm}\label{thm8} {\rm (i)} Assume that $ \max\{2^*(\alpha_2),2_*(\beta)\} <p<2^*(\beta)$ with $\beta>-2$. Then for any $\mu<0$ being fixed, there exists a $\lambda^{**}>0$ such that for   $\lambda>\lambda^{**}$, \eqref{4.1} has  a  nonnegative ground state solution.

 {\rm (ii)}  Assume that $p$ satisfies
 \begin{equation}\label{4.16}
 \left\{
 \begin{array}{ll}
 (1) \ \max\{2(3+\alpha_2),2(2+\beta)\}<p<2^*(\beta), & \beta \geqslant -1, N=3, \\
 (2) \ 6+2\alpha_2<p<2^*(\beta), & \beta< -1, N=3, \\
(3) \ 2^*(\alpha_2)<p<2^*(\beta),&  -2<\beta\leqslant 0, N\geqslant4,  \\
(4) \ \max\{2^*(\alpha_2),\frac{2(N-2+\beta)}{N-2}\}<p<2^*(\beta), &  \beta>0, \ N\geqslant4.
\end{array}
\right.
\end{equation}
Then for any $\lambda>0$, there exists a $\mu^{*}<0$ such that for any $\mu^{*}<\mu<0$, the problem \eqref{4.1} possesses a nonnegative ground state solution.

(iii) Assume that $p$ satisfies \eqref{new3.6}. Then for any $0<\lambda<\lambda_{1\beta}$, there exists $\mu^{**}<0$ such that for any $\mu^{**}<\mu<0$, the problem \eqref{4.1} possesses a nonnegative ground state solution.
\end{thm}

We will work with the functional $\Phi_{\lambda, \mu}$ defined by \eqref{4.2} and use the corresponding notations given by \eqref{4.3} and \eqref{4.4}. For convenience, we give a assumption
\begin{eqnarray}
\left\{
  \begin{array}{ll}
\max\{2^*(\alpha_2),2_*(\beta)\} <p<2^*(\beta)\ {\rm and}\ \beta>-2, \lambda>0;\\
p=2\ {\rm and}\ -2<\beta<0, 0<\lambda<\lambda_{1\beta}.
  \end{array}
\right.\label{4.17}
\end{eqnarray}
 We first have
\begin{pro}\label{pro4.8} \ Let \eqref{4.17} hold.
\begin{enumerate}
 \item[{\rm (i) }] For $u\in H_r^1(\mathbb{R}^N)\backslash\{0\}$, there exists a unique $t_u>0$ such that  $t_u u\in\mathcal{N}_{\lambda,\mu} $ and $\Phi_{\lambda,\mu} (t_u u)=\max_{t\geqslant  0}\Phi_{\lambda,\mu} (tu)$, moreover, the manifold $\mathcal{N}_{\lambda,\mu}$ is nonempty;
     \item[{\rm (ii) }] The manifold $\mathcal{N}_{\lambda,\mu}$ is $C^1$ regular, that is $\langle \Psi_{\lambda,\mu}'(u), u\rangle=\langle \Phi_{\lambda,\mu}''(u)u, u\rangle\ne0$ for any $u\in \mathcal{N}_{\lambda,\mu}$;
     \item[{\rm (iii) }] $\mathcal{N}_{\lambda,\mu}$ is closed and bounded away from $0$, and $m_{\lambda,\mu}>0$.
 \end{enumerate}
 \end{pro}
 It follows also the arguments of \cite[Chapter 4]{1996Willem} that
 \begin{eqnarray}  \label{4.18}
 m_{\lambda,\mu} = c_{\lambda,\mu}:=\inf_{u\in H_r^1(\mathbb{R}^N) \setminus\{0\}}
  \max_{t\geqslant 0} \Phi_{\lambda,\mu} (t u).
 \end{eqnarray}

 For  the sake of convenience, we set
 \begin{equation*} \begin{array}{ll} A(u):=\|u\|^2, & B(u):=\|u\|^{2^*(\alpha_1)}_{L^{2^*(\alpha_1)}(\mathbb{R}^N, |x|^{\alpha_1})}, \\
   C(u) :=\|u\|^{2^*(\alpha_2)}_{L^{2^*(\alpha_2)}(\mathbb{R}^N, |x|^{\alpha_2})}, &  D(u):=\|u\|^p_{L^p(\mathbb{R}^N, |x|^{\beta})}.
 \end{array} \end{equation*} Then $\Phi_{\lambda,\mu}$ can be rewritten as
 $$ \Phi_{\lambda,\mu}(u)= \frac12\|u\|^2-\frac{1}{2^*(\alpha_1)} B(u)- \frac{\mu}{2^*(\alpha_2)}C(u)-\frac{\lambda}{p} D(u).$$

\begin{lem} \label{lem5} Assume \eqref{4.17} holds.
Any a $(PS)_c$ sequence $\{u_n\}\subset H^1_r(\mathbb{R}^N)$ of $\Phi_{\lambda,\mu}$ with
 $$c \leqslant c^*:=\frac{2+\alpha_1}{2(N+\alpha_1)}S_{\alpha_1}^{\frac{N+\alpha_1}{2+\alpha_1}}$$
has  a convergent subsequence.
\end{lem}
{\bf Proof.} \ For any a ${\rm(PS)_c}$ sequence $\{u_n\}\subset H^1_r({\mathbb{R}^N})$, we have that as $n\rightarrow\infty$,
 \begin{eqnarray}
\Phi_{\lambda, \mu}(u_n)&=&\frac{1}{2}\|u_n\|^2-\frac{1}{2^*(\alpha_1)}B(u_n)
-\frac{\mu}{2^*(\alpha_2)}C(u_n)-\frac{\lambda}{p}D(u_n) \to c,\label{4.19}\\
\langle \Phi_{\lambda, \mu}'(u_n), u_n\rangle &= &\|u_n\|^2-B(u_n)-C(u_n)-\lambda D(u_n)=o(1)\|u_n\|.\label{4.20}
 \end{eqnarray}
 As $p>2^*(\alpha_2)>2$, we can deduce from \eqref{4.19}, \eqref{4.20} and combining with \eqref{3.11} that
 $\{u_n\}$ is bounded in $H^1_r(\mathbb{R}^N)$.  Up to a subsequence if necessary,  we can assume that there is $u\in  H^1_r({\mathbb{R}^N})$ such that $u$ satisfies
 \begin{eqnarray}
  \left\{\begin{array}{ll}
&u_n\rightharpoonup u \ \ \mbox{in} \ H^1_r(\mathbb{R}^N);\\
&u_n\rightharpoonup u \ \ \mbox{in} \ L^{2^*(\gamma)}(\mathbb{R}^N;|x|^\gamma), \ \gamma>-2;\\
&u_n\rightarrow u \ \  \mbox{in} \ L^p(\mathbb{R}^N;|x|^\beta),  \ \   \ 2_*(\beta)<p<2^*(\beta),  \ \beta>-2;\\
&u_n\rightarrow u \ \  \mbox{in} \ L^2(\mathbb{R}^N;|x|^\beta),  \ \ \mbox{for} \ p=2_*(\beta)=2 \ \hbox{with} \ -2<\beta<0;\\
&u_n(x)\rightarrow u(x) \ \ \mbox{a.e.\ on} \  \mathbb{R}^N.
  \end{array}
\right.\label{new4.21}
\end{eqnarray}
 Furthermore, we can get that $u$ satisfies weakly the equation
 $$-\Delta u+u=|x|^{\alpha_1}|u|^{2^*(\alpha_1)-2}u+\mu|x|^{\alpha_2}|u|^{2^*(\alpha_2)-2}u+\lambda |x|^\beta|u|^{p-2}u$$ so that
 \begin{eqnarray} \label{4.21}
 \langle \Phi'_{\lambda,\mu}(u), u\rangle=0, \ \ \ \  \ \
 \Phi_{\lambda,\mu}(u) \geqslant0.
 \end{eqnarray}
 Set $v_n:=u_n-u$. By  \eqref{new4.21},\eqref{4.20} and Brezis-Lieb lemma (\cite{1983B-Lieb}), we have
\begin{equation}\label{4.22}
\|v_n\|^2-B(v_n)-\mu C(v_n)\rightarrow0 \ \ \ n\to\infty.
\end{equation}
Up to a subsequence if necessary,  we assume that
$$\displaystyle\lim_{n\rightarrow\infty}A(v_n)=A_\infty, \ \ \ \lim_{n\rightarrow\infty} B(v_n)=B_\infty, \ \ \ \lim_{n\rightarrow\infty}
C(v_n)=C_\infty.$$
By \eqref{4.22}, we have
$$A_\infty-\mu C_\infty=B_\infty.$$
Assume that $A_\infty>0$. Then we can prove by using  Lemma \ref{lem4.1} and Lemma \ref{lem4.2} that $B_\infty>0 $ and $C_\infty>0$.
It follows from H\'enon-Sobolev inequality and $\mu<0$ that
\begin{eqnarray*}
S_{\alpha_1} B_\infty^{\frac{2}{2^*(\alpha_1)}}\leqslant A_\infty \leqslant A_\infty-\mu C_\infty=B_\infty.
\end{eqnarray*}
Thus
$$B_\infty \geqslant S_{\alpha_1}^{\frac{2^*(\alpha_1)}{2^*(\alpha_1)-2}}, \ \ \ \ A_\infty\geqslant  S_{\alpha_1}^{\frac{N+\alpha_1}{2+\alpha_1}}.$$
By \eqref{4.19} and \eqref{new4.21} we deduce that
\begin{eqnarray*}
  \Phi_{\lambda,\mu}(u)+\frac{1}{2}A(v_n) -\frac{1}{2^*(\alpha_1)}B(v_n)-\frac{\mu}{2^*(\alpha_2)}C(v_n) \rightarrow c \ \ \ n\rightarrow\infty.
\end{eqnarray*}
Therefore we have by \eqref{4.21} that
\begin{eqnarray*}
c \geqslant \frac{A_\infty}{2}-\frac{B_\infty}{2^*(\alpha_1)}-\frac{\mu C_\infty}{2^*(\alpha_2)}>\frac{2+\alpha_1}{2(N+\alpha_1)}S_{\alpha_1}^{\frac{N+\alpha_1}{2+\alpha_1}}.
\end{eqnarray*}
which contradicts the choice of $c$. Hence $A_\infty=0$ and the proof is complete. \hfill $\Box$

  \begin{lem}\label{lem8} \ {\rm (i)} Assume that $ \max\{2^*(\alpha_2),2_*(\beta)\} <p<2^*(\beta)$ with $\beta>-2$. Then for any $\mu<0$ being fixed, there exists a $\lambda^{**}>0$ such that for   $\lambda>\lambda^{**}$, \begin{eqnarray}
  m_{\lambda, \mu} \leqslant \frac{2+\alpha_1}{2(N+\alpha_1)}S_{\alpha_1}^{\frac{N+\alpha_1}{2+\alpha_1}}. \label{4.23}
  \end{eqnarray}
  {\rm (ii)} Assume that $p$ satisfies \eqref{4.16}.  Then for any $\lambda>0$ being fixed,
  there exists  $\mu^*<0$ such that for any $\mu^*<\mu<0$ \eqref{4.23} holds.\\
{\rm(iii)} Assume that \eqref{new3.6} holds. Then for any $0<\lambda<\lambda_{1\beta}$, there exists $\mu^{**}<0$ such that for $\mu^{**}<\mu<0$, \eqref{4.23} holds.
  \end{lem}
{\bf Proof} \ By \eqref{4.18}, we only need to find a nonzero $v \in H_r^1(\mathbb{R}^N)$ such that
  \begin{eqnarray} \max_{t \geqslant0} \Phi_{\lambda, \mu}(t v )
   \leqslant \frac{2+\alpha_1}{2(N+\alpha_1)}S_{\alpha_1}^{\frac{N+\alpha_1}{2+\alpha_1}}. \label{4.24}
  \end{eqnarray}
The proof of the case (i) is similar to that of Lemma \ref{lem4.6}. We prove the case (ii), (iii) by using the function
  \begin{eqnarray*}
u_{\epsilon,\alpha_1}= \varphi(|x|) U_{\epsilon,\alpha_1}
  \end{eqnarray*}
  defined by \eqref{2.7} with $\alpha$ being replaced by $\alpha_1$. We have the following estimates:
 \begin{eqnarray}
 &&\int_{\mathbb{R}^N}\left|\nabla u_{\epsilon,{\alpha_1}}\right|^2dx
 =S_{\alpha_1}^{\frac{N+{\alpha_1}}{2+{\alpha_1}}}+O\left(\epsilon^{N-2}\right).\label{4.25}\\
 &&\int_{\mathbb{R}^N}|x|^{\alpha_1} |u_{\epsilon,{\alpha_1}}|^{2^*({\alpha_1})}dx
 =S_{\alpha_1}^{\frac{N+{\alpha_1}}{2+{\alpha_1}}}+O\left(\epsilon^{N+{\alpha_1}}\right).\label{4.26}
 \end{eqnarray}
 \begin{eqnarray}
 \int_{\mathbb{R}^N}|u_{\epsilon,{\alpha_1}}|^2dx=\left \{\begin{array}{ll}
 O(\epsilon), &  N=3;\\
 O(\epsilon^2)+C\epsilon^2|\ln \epsilon|, & N=4;\\
 O(\epsilon^{N-2})+C\epsilon^2,& N \geqslant 5.
 \end{array}
 \right.\label{4.27}
\end{eqnarray}
\begin{eqnarray}
\int_{\mathbb{R}^N}|x|^{\alpha_2}|u_{\epsilon,\alpha_1}|^{2^*(\alpha_2)}dx
=\widetilde{K}+O(\epsilon^{N+\alpha_2}),\label{4.28}
\end{eqnarray}
where $\widetilde{K}=\int_{\mathbb{R}^N} |x|^{\alpha_2}|U_{1,\alpha_1}|^{2^*(\alpha_2)}dx$.
\begin{eqnarray}
\int_{\mathbb{R}^N}|x|^\beta|u_{\varepsilon,{\alpha_1}}|^pdx=\left \{\begin{array}{ll}
O(\epsilon^{\frac{N+\beta}{2}})+C\epsilon^{\frac{N+\beta}{2}}|\ln \epsilon|, & p=\frac{N+\beta}{N-2};\\
 O(\epsilon^{\frac{p(N-2)}{2}})+C\epsilon^{N+\beta-\frac{p(N-2)}{2}}, & p>\frac{N+\beta}{N-2}.
\end{array}
\right.\label{4.29}
\end{eqnarray}
For the case (ii). If  $p$ satisfies \eqref{4.16},  we have by \eqref{4.25}--\eqref{4.29} that
\begin{eqnarray}\label{4.30}
 \begin{array}{ll}
 \displaystyle \lim_{\epsilon\rightarrow ^+} A(u_{\epsilon,\alpha_1})=S_{\alpha_1}^{\frac{N+{\alpha_1}}{2+{\alpha_1}}}, \ \ \  \ &  \displaystyle \lim_{\epsilon\rightarrow 0^+} B(u_{\epsilon,\alpha_1})=S_{\alpha_1}^{\frac{N+{\alpha_1}}{2+{\alpha_1}}},\\
 \displaystyle \lim_{\epsilon\rightarrow 0^+}C(u_{\epsilon,\alpha_1})=\widetilde{K},  \ \ \ \   & \displaystyle \lim_{\epsilon\rightarrow 0^+} D(u_{\epsilon,\alpha_1})=0.
 \end{array}
 \end{eqnarray}
Let $\lambda>0$ be fixed. Then there is a unique  $t_{\epsilon,\mu}>0$ such that
 \begin{eqnarray}
\displaystyle && \sup_{t\geqslant0} \Phi_{\lambda, \mu}(t u_{\epsilon,\alpha_1}) = \Phi_{\lambda, \mu}(t_{\epsilon,\mu} u_{\epsilon,\alpha_1}) \nonumber\\
&=&\frac{t_{\epsilon,\mu}^2}{2}A(u_{\epsilon,\alpha_1})
-\frac{t_{\epsilon,\mu}^{2^*(\alpha_1)}}{2^*(\alpha_1)}B(u_{\epsilon,\alpha_1})
-\frac{\mu t_{\epsilon,\mu}^{2^*(\alpha_2)}}{2^*(\alpha_2)}
C(u_{\epsilon,\alpha_1})
-\frac{\lambda t_{\epsilon,\mu}^p}{p}D(u_{\epsilon,\alpha_1}).\label{4.31}
 \end{eqnarray}
Moreover, $t_{\epsilon,\mu} u_{\epsilon,\alpha_1} \in \mathcal{N}_{{\lambda,\mu}}$ and
\begin{eqnarray}
A(u_{\epsilon,\alpha_1})
-t_{\epsilon,\mu}^{2^*(\alpha_1)-2}B(u_{\epsilon,\alpha_1})
-\mu t_{\epsilon,\mu}^{2^*(\alpha_2)-2}
C(u_{\epsilon,\alpha_1})
-\lambda t_{\epsilon,\mu}^{p-2}D(u_{\epsilon,\alpha_1})=0.\label{4.32}
\end{eqnarray}
For $p \geqslant \frac{N+\beta}{N-2}$, by \eqref{4.30} we have  for $\epsilon>0$ small enough that
\begin{eqnarray*}
A(u_{\epsilon,\alpha_1})\leqslant  2S_{\alpha_1}^{\frac{N+\alpha_1}{2+\alpha_1}}, \ \ B(u_{\epsilon,\alpha_1}) \geqslant  \frac{1}{2}S_{\alpha_1}^{\frac{N+\alpha_1}{2+\alpha_1}}, \ \
C(u_{\epsilon,\alpha_1})\leqslant 2 {\widetilde{K}}.
\end{eqnarray*}
For $-T< \mu <0 $ with $T>0$ and for $\epsilon>0$ small enough, by \eqref{4.32}, we have
\begin{eqnarray}
0& \leqslant & A(u_{\epsilon,\alpha_1})-t_{\epsilon,\mu}^{2^*(\alpha_1)-2}B(u_{\epsilon,\alpha_1})
+T t_{\epsilon,\mu}^{2^*(\alpha_2)-2}C(u_{\epsilon,\alpha_1})\nonumber\\
& \leqslant & 2 S_{\alpha_1}^{\frac{N+\alpha_1}{2+\alpha_1}}-\frac{1}{2}t_{\epsilon,\mu}^{2^*(\alpha_1)-2}
S_{\alpha_1}^{\frac{N+\alpha_1}{2+\alpha_1}}+
2 T t_{\epsilon,\mu}^{2^*(\alpha_2)-2} \tilde{K}:=g(t_{\epsilon,\mu}).\label{4.33}
\end{eqnarray}
It is easy to see that there exists $s_0>0$ such that $g(s)$ is increasing on $(0,s_0)$ and $g(s)$ is decreasing on $(s_0,\infty)$. By  \eqref{4.33}, $g(0)=2S_{\alpha_1}^{\frac{N+\alpha_1}{2+\alpha_1}}>0$ and $\displaystyle
\lim_{s\rightarrow \infty}g(s)=-\infty$.
 It follows from \eqref{4.30} and \eqref{4.32} that there exists $M>0$ independent on $\epsilon>0$ and $\mu<0$ such that
\begin{eqnarray}
 0< t_{\epsilon, \mu}\leqslant  M,  \ \ \ {\rm for} \ \epsilon>0, \ -\mu>0 \ \  {\rm small\ enough}.\label{4.34}
\end{eqnarray}
Furthermore,  by \eqref{4.32} and \eqref{4.34}, we get
\begin{eqnarray}
 t_{\epsilon,\mu}\rightarrow 1 \ \ \     \epsilon\rightarrow0^+, \ \mu\rightarrow0^-.\label{4.35}
\end{eqnarray}
It follows from \eqref{4.31}, we have
\begin{eqnarray*}
\displaystyle \sup_{t\geqslant 0} \Phi_{\lambda,\mu}(t u_{\epsilon,\alpha_1}) \leqslant \max_{t\geqslant0}\left\{\frac{t^2}{2}A(u_{\epsilon,\alpha_1})
-\frac{t^{2^*(\alpha_1)}}{2^*(\alpha_1)}B(u_{\epsilon,\alpha_1})-\frac{\lambda t^p}{p}D(u_{\epsilon,\alpha_1})\right\}
-\frac{\mu t_{\epsilon,\mu}^{2^*(\alpha_2)}}{2^*(\alpha_2)}C(u_{\epsilon,\alpha_1}).
\end{eqnarray*}
By the proof of Lemma \ref{lem2}, we have for $p$ satisfying \eqref{4.16} and for $\epsilon>0$ small enough,
\begin{eqnarray}
\max_{t\geqslant 0}\left\{\frac{t^2}{2}A(u_{\epsilon,\alpha_1})
-\frac{t^{2^*(\alpha_1)}}{2^*(\alpha_1)}B(u_{\epsilon,\alpha_1})-\frac{\lambda t^p}{p}D(u_{\epsilon,\alpha_1})\right\}
<\frac{2+\alpha_1}{2(N+\alpha_1)}S_{\alpha_1}^{\frac{N+\alpha_1}{2+\alpha_1}}.\label{4.36}
\end{eqnarray}
Since from \eqref{4.30} and \eqref{4.35} we can deduce that
\begin{eqnarray}
\frac{\mu t_{\epsilon,\mu}^{2^*(\alpha_2)}}{2^*(\alpha_2)}C(u_{\epsilon,\alpha_1})\rightarrow0, \ \ \epsilon \to0^+, \ \ \mu\to0^-,\label{4.37}
\end{eqnarray}
 there exists $\mu^*<0$ such that $u_{\epsilon, \alpha_1}$ satisfies \eqref{4.24} for any $\mu^*<\mu<0$.

For the case of (iii). The estimation \eqref{4.29} with $p=2$ becomes
\begin{eqnarray}
\int_{\mathbb{R}^N}|x|^\beta|u_{\varepsilon,{\alpha_1}}|^2dx=\left \{\begin{array}{ll}
O(\epsilon^{\frac{N+\beta}{2}})+C\epsilon^{\frac{N+\beta}{2}}|\ln \epsilon|, & \beta=N-4;\\
 O(\epsilon^{\frac{p(N-2)}{2}})+C\epsilon^{N+\beta-\frac{p(N-2)}{2}}, & -2<\beta<N-4.
\end{array}
\right.\label{4.38}
\end{eqnarray}
Using the similar arguments of (ii), combining with the inequality \eqref{3.11} and Lemma \ref{lem2}, we know that \eqref{4.35} and \eqref{4.36} hold. The fact \eqref{4.37} implies there exists $\mu^{**}<0$ such that \eqref{4.24} holds for any $\mu^{**}<\mu<0$. The proof  is complete. \hfill$\Box$

\noindent{\bf Proof of Theorem \ref{thm8}} \ The argument is similar to that of Theorem \ref{thm6} and Theorem \ref{thm5}. We omit the details.  $\hfill\Box$

We finish this paper by pointing out that one may consider the equations with much more critical exponents. We leave   the precise statements for the interested readers. In a forthcoming paper we will consider the critical H\'enon-Sobolev exponent problems under the perturbations of  general functions.


\begin{thebibliography}{99}

\footnotesize

 \bibitem{1973AR}   A. Ambrosetti and P. H. Rabinowitz, Dual variational methods  in  critical point theory and applications. \emph{  J. Funct. Anal.,}
    \textbf{14}(1973), 349--381.

\bibitem{1983B-Lions} H. Berestycki, P-L. Lions, Nonlinear scalar field equations. I. Existence of a ground state. \emph{Arch. Rational Mech. Anal.,} \textbf{82} (1983) 313--345.

\bibitem{1983B-Lieb}  H. Br\'{e}zis, E. Lieb, A relation between pointwise convergence of functions and convergence of functionals. \emph{Proc. Amer. Math. Soc.,} \textbf{88} (1983) 486--490.

\bibitem{1983BN}  H. Br\'{e}zis, L.  Nirenberg, Positive solutions of nonlinear elliptic equations involving critical Sobolev exponents. \emph{Comm. Pure Appl. Math.,} \textbf{36} (1983) 437--477.

\bibitem{2001CW}  F. Catrina, Z.-Q. Wang, On the Caffarelli-Kohn-Nirenberg inequalities: sharp constants, existence (and nonexistence), and symmetry of extremal functions. \emph{Comm. Pure Appl. Math.,} \textbf{54} (2001) 229--258.


\bibitem{2010CL} J.-L. Chern, C. -S. Lin, Minimizers of Caffarelli-Kohn-Nirenberg inequalities with the singularity on the boundary. \emph{Arch. Ration. Mech. Anal.,} \textbf{197} (2010) 401--432.

\bibitem{1992Egnell} H. Egnell, Positive solutions of semilinear equations in cones. {\em Trans. Amer. Math. Soc.,} \textbf{330}(1992) 191--201.

    \bibitem{2001FG} A. Ferrero, F. Gazzola, Existence of solutions for singular critical growth semilinear elliptic equations. \emph{J. Differential Equations,} \textbf{177} (2001)  494--522.

     \bibitem{2004GK} N. Ghoussoub, X. Kang, Hardy-Sobolev critical elliptic equations with boundary singularities. \emph{Ann. Inst. H. Poincar\'e Anal. Non Lin\'eaire,} \textbf{21} (2004) 767--793.


\bibitem{2016GR1}  N. Ghoussoub, F. Robert, Sobolev inequalities for the Hardy-Schr\"odinger operator:
 extremals and critical dimensions. {\em  Bull. Math. Sci.,} {\bf 6}(2016) 89--144.



\bibitem{2017GR2} N.Ghoussoub, F. Robert, Hardy-singular boundary mass and Sobolev-critical variational problems. {\em Anal. PDE,} \textbf{10}(2017) 1017--1079.

 \bibitem{2000GY} N. Ghoussoub, C. Yuan, Multiple solutions for quasi-linear PDEs involving the critical Sobolev and Hardy exponents. {\em Trans. Amer. Math. Soc.,} \textbf{352}(2000) 5703--5743.



 \bibitem{1981GS} B. Gidas, J. Spruck, Global and local behavior of positive solutions of nonlinear elliptic equations. \emph{Comm. Pure Appl. Math.,} \textbf{24} (1981) 525--598.


\bibitem{2013GGN} F. Gladiali, M. Grossi, S. L. N. Neves, Nonradial solutions for the H\'{e}non equation in $\mathbb{R}^N$. \emph{Adv. Math.,} \textbf{249} (2013) 1--36.


    \bibitem{1973Henon} M. H\'{e}non, Numerical experiments on the stability oh spherical stellar systems. \emph{Astronom. Astrophys.,} \textbf{24} (1973), 229--238.


\bibitem{2010HL} C.-H. Hsia, C.-S. Lin, H. Wadade, Revisiting an idea of Br\'{e}zis and Nirenberg. \emph{J. Funct. Anal.,} \textbf{259} (2010) 1816--1849.


 \bibitem{2011HL} C.-S. Lin, C.H. Hsia, H. Wadade, The existence of positive solutions to the semilinear elliptic equation involving the Sobolev and the Sobolev-Hardy critical terms. Harmonic analysis and nonlinear partial differential equations, 133--157, {\em RIMS K\^{o}ky\^{u}roku Bessatsu, B26,} Res. Inst. Math. Sci. (RIMS), Kyoto, 2011.

 \bibitem{2012Li-Lin} Y. Y. Li, C.-S. Lin, A nonlinear elliptic PDE and two Sobolev-Hardy critical exponents. {\em Arch. Ration. Mech. Anal.,} \textbf{203}(2012) 943--968.



 \bibitem{1983Lieb} E. Lieb, Sharp constants in the Hardy--Littlewood--Sobolev and related inequalities.
 \emph{Ann. of Math.} \textbf{118} (1983) 349--374.







\bibitem{1982Ni}  W. M. Ni, A nonlinear Dirichlet problem on the unit ball and its applications. \emph{Indiana Univ. Math. J.,} \textbf{31} (1982) 801--807.

\bibitem{1965Pohozaev} S. I. Poho\u{z}aev, On the eigenfunctions of the equation $\Delta u+\lambda f(u)=0$.(Russian) \emph{Dokl. Akad. Nauk,} \textbf{165} (1965) 36--39.





 \bibitem{2002SSW} D. Smets, J. Su, M. Willem,   Non-radial ground states for the H\'{e}non equation. \emph{ Commun. Contemp. Math.,} \textbf{4} (2002) 467--480.



    \bibitem{1977Strauss}  W.A. Strauss, Existence of solitary waves in higher dimensions. \emph{Comm. Math. Phys.,} \textbf{55}(1977) 149--162.




\bibitem{2008Struwe} M. Struwe, {\it Variational methods. Applications to nonlinear partial differential equations and Hamiltonian systems,} Fourth edition, Springer-Verlag, Berlin, 2008.



 \bibitem{2007SWW-1} J.  Su, Z.-Q. Wang, M. Willem, Nonlinear Schr\"{o}dinger equations with unbounded and decaying radial potentials. \emph{Commun. Contemp. Math.,} \textbf{9} (2007) 571--583.

\bibitem{2007SWW-2}  J.  Su, Z.-Q. Wang, M. Willem, Weighted Sobolev embedding with unbounded and decaying radial potential. \emph{J. Differential Equations,} \textbf{238} (2007) 201--219.

\bibitem{2011SW} J. Su, Z.-Q. Wang, Sobolev type embedding and quasilinear elliptic equations with radial potentials. \emph{J. Differential Equations,} \textbf{250}(2011) 223--242.

  \bibitem{1976Talenti} G. Talenti, Best constant in Sobolev inequality. {\em Ann. Mat. Pura Appl.,} \textbf{110}(1976) 353--372.



    \bibitem{2019Wang-Su} C. Wang, J.  Su, Positive radial solutions of critical H\'{e}non equations on the unit ball in $\mathbb{R}^N$. \emph{ Math. Methods Appl. Sci.,} \textbf{45}(2022) 11769--11806.




\bibitem{1996Willem} M. Willem, {\em Minimax Theorems}. {Birkh\"auser Boston. Inc. Boston,} 1996.

















\end{thebibliography}
\end{document}